\documentclass{article}
\usepackage{color}
\usepackage{amsfonts,amsmath,upref,amssymb,amsthm,amsxtra,amscd,enumerate}
\catcode`\@=11 \@addtoreset{equation}{section}

\catcode`\@=12

\newtheorem{Theorem}{Theorem}[section]
\newtheorem{Proposition}{Proposition}[section]
\newtheorem{Lemma}{Lemma}[section]
\newtheorem{Corollary}{Corollary}[section]
\newtheorem{Definition}{Definition}[section]
\newtheorem{Remark}{Remark}[section]

\newcommand{\bTheorem}[1]{
\begin{Theorem} \label{T#1} }
\newcommand{\eT}{\end{Theorem}}

\newcommand{\bProposition}[1]{
\begin{Proposition} \label{P#1}}
\newcommand{\eP}{\end{Proposition}}

\newcommand{\bLemma}[1]{
\begin{Lemma} \label{L#1} }
\newcommand{\eL}{\end{Lemma}}

\newcommand{\bCorollary}[1]{
\begin{Corollary} \label{C#1} }
\newcommand{\eC}{\end{Corollary}}

\newcommand{\bFormula}[1]{
\begin{equation} \label{#1}}
\newcommand{\eF}{\end{equation}}

\newcommand{\Ov}[1]{\overline{#1}}
\newcommand{\DC}{C^\infty_c}

\newcommand{\vr}{\varrho}
\newcommand{\vre}{\vr_\ep}
\newcommand{\vte}{\vt_\ep}
\newcommand{\vue}{\vu_\ep}

\newcommand{\vt}{\vartheta}
\newcommand{\vu}{\vc{u}}
\newcommand{\vc}[1]{{\bf #1}}
\newcommand{\Div}{{\rm div}_x}
\newcommand{\Grad}{\nabla_x}

\newcommand{\tn}[1]{\mbox {\F #1}}
\newcommand{\dx}{\,{\rm d} {x}}
\newcommand{\vph}{{\boldsymbol \varphi}}
\newcommand{\vq}{\mathbf{q}}
\newcommand{\dt}{\,{\rm d} t }

\newcommand{\dxdt}{\,{\rm d}x{\rm d}t}

\newcommand{\K}{{\cal K}}
\newcommand{\ep}{\varepsilon}

\font\F=msbm10 scaled 1000
\newcommand{\R}{\mbox{\F R}}

\providecommand{\abs}[1]{\left\lvert#1\right\rvert}
\newcommand{\sil}{\rightarrow}

\definecolor{grey}{rgb}{0.85,0.85,0.85}
\long\def\greybox#1{%
    \newbox\contentbox%
    \newbox\bkgdbox%
    \setbox\contentbox\hbox to \hsize{%
        \vtop{
            \kern\columnsep
            \hbox to \hsize{%
                \kern\columnsep%
                \advance\hsize by -2\columnsep%
                \setlength{\textwidth}{\hsize}%
                \vbox{
                    \parskip=\baselineskip
                    \parindent=0bp
                    #1
                }%
                \kern\columnsep%
            }%
            \kern\columnsep%
        }%
    }%
    \setbox\bkgdbox\vbox{
        \color{grey}
        \hrule width  \wd\contentbox %
               height \ht\contentbox %
               depth  \dp\contentbox
        \color{black}
    }%
    \wd\bkgdbox=0bp%
    \vbox{\hbox to \hsize{\box\bkgdbox\box\contentbox}}%
    \vskip\baselineskip%
}

\begin{document}

%%%%%%%%%%%%%%%%%%%%%%%%%%%%%%%%

\title{Flow of heat conducting fluid in a time dependent domain}

\author{
Ond{\v r}ej Kreml $^{1}$\thanks{The work of O.K.and  \v S.N. was supported by 7AMB16PL060 and by RVO 67985840.  Stay of V. M. in Imperial College London was supported by GA17-01747S and RVO: 67985840. Stay of O.K. at Imperial College London was supported by the grant Iuventus Plus 0871/IP3/2016/74.}
\and V\'aclav M\'acha $^{1}$\footnotemark[1]
\and {\v S}{\' a}rka Ne{\v c}asov{\' a} $^1$\footnotemark[1]
\and Aneta Wr\'oblewska-Kami\'nska $^2$\thanks{The work of A.W.-K.
is partially supported by a Newton Fellowship of the
Royal Society and by the grant Iuventus Plus 0871/IP3/2016/74 of Ministry of Sciences and Higher
Education RP. Her stay at Institute of Mathematics of Academy of Sciences, Prague was supported by 7AMB16PL060.}
}

\maketitle

\bigskip

\centerline{$^1$Institute of Mathematics of the Academy of Sciences of
the Czech Republic} \centerline{\v Zitn\' a 25, 115 67 Praha 1,
Czech Republic}
\bigskip

\centerline{$^2$  Institute of Mathematics, Polish Academy of Sciences}
\centerline{\'Sniadeckich 8, 00-656 Warszawa, Poland}
\centerline{and Department of Mathematics, Imperial College London}
\centerline{London SW7 2AZ, United Kingdom}

\medskip

\begin{abstract}

We consider a flow of heat conducting fluid inside a moving domain whose shape in time is prescribed. The flow in this case is governed by the Navier-Stokes-Fourier system consisting of equation of continuity, momentum balance, entropy balance and energy equality. The velocity is supposed to fulfill the full-slip boundary condition and we assume that the fluid is thermally isolated. In the presented article we show the existence of a variational solution.
\end{abstract}

\medskip


\medskip

\section{Introduction}\label{i}

The flow of heat conducting fluid inside a moving domain is an interesting problem with a lot of practical applications and deserves its attention. Let us mention modeling of the motion of a piston in a cylinder filled by a viscous heat conducting gas. There are many references on this problem in the statistical physics. We can mention works of Lieb \cite{Lieb}, Gruber et al. \cite{Gruber1,Gruber2,Gruber3}, Wright \cite {Wright1, Wright2, Wright3}, etc. The problem was investigated by Shelukhin \cite{Shel_first}, Antman and Wilber \cite{AW} in the case of homogeneous boundary conditions for barotropic case. The extension to the case of non-homogeneous boundary conditions can be found in the work of  Maity et al. \cite{DTT16}. The motion of a piston in a cylinder filled by a viscous heat conducting gas was studied by Shelukhin \cite{Shel}. His results coincide with the statistical mechanics scenario for the thermally insulting piston. Rather complex behavior of piston was proven by Feireisl et al. \cite{fe}.

Although the full Navier-Stokes-Fourier system in a moving domain has been already investigated (see \cite{KrMaNeWr}) it turns out that it is more useful to examine the following form of the NSF system (under more general hypotheses on the pressure) at least from the point of view of further analysis like asymptotic limits, dimension reduction and so on:
\bFormula{i1a}
\partial_t \vr + \Div (\vr \vu) = 0,
\eF
\bFormula{i1b}
\partial_t (\vr \vu) + \Div (\vr \vu \otimes \vu) + \Grad p(\vr,\vt) = \Div \tn{S}(\Grad \vu),
\eF
\bFormula{i1c}
\partial_t (\vr s) + \Div (\vr s \vu) + \Div \left(\frac{\vq}{\vt}\right) = \sigma,
\eF
\bFormula{i1d}
\frac{{\rm d}}{{\rm d}t} \int \left( \vr |\vu|^2 + \vr e \right)\, \dx = 0.
\eF
These equations, which are considered on a time-space domain $(0,T)\times \Omega_t\subset (0,\infty)\times \mathbb R^3$ where $\Omega_t$ is a time dependent domain, are mathematical formulations of the balance of mass, linear momentum, entropy and total energy respectively. Unknowns are the density $\vr:(0,T)\times \Omega_t\mapsto [0,\infty)$, the velocity $\vu:(0,T)\times \Omega_t \mapsto \mathbb R^3$ and the temperature $\vt:(0,T)\times \Omega_t\mapsto [0,\infty)$. Other quantities appearing in these equations are functions of the unknowns, namely the stress tensor $\tn{S}$, the internal energy $e$, the pressure $p$, the entropy $s$, and the entropy production rate $\sigma$. Their needed properties are mentioned later on. For simplicity we do not consider in this work any external forces.

The time dependent domain $\Omega_t$ is prescribed by movement of its boundary on the time interval $[0,T]$. Namely, the boundary of the domain $\Omega_t$ occupied by the fluid is described by a given velocity field $\vc{V}(t,x)$ where $t \geq 0$ and $x \in \mathbb{R}^3$. Supposing $\vc{V}$ is regular enough we can associate the following system of equations
\begin{equation*}
\frac{{\rm d}}{{\rm d}t} \vc{X}(t, x) = \vc{V} \Big( t, \vc{X}(t, x) \Big),\ t > 0,\ \vc{X}(0, x) = x,
\end{equation*}
with our domain, then we set
\[
\Omega_\tau = \vc{X} \left( \tau, \Omega_0 \right), \ \mbox{where} \ \Omega_0 \subset \mathbb{R}^3 \ \mbox{is a given domain},\
\Gamma_\tau = \partial \Omega_\tau, \ \mbox{and}\ Q_\tau = \cup_{t\in(0,\tau)} \{t\}  \times \Omega_t .%\{ (t,x) \ |\ t \in (0,\tau), \ x \in \Omega_t \}.
\]
We assume that the volume of the domain can not degenerate in time, namely
\begin{equation}\label{V_0}
\mbox{ there exists }V_0 >0 \mbox{ such that } |\Omega_\tau| \geq V_0 \mbox{ for all } \tau \in [0,T].
\end{equation}
Moreover, we assume that
\begin{equation}\label{Vselenoidal}
\Div \vc{V} = 0\ \mbox{on the neighborhood of }\Gamma_\tau,
\end{equation}
see Remark \ref{divboundry}.

We assume that the boundary of the physical domain is impermeable. This is described by the condition
\bFormula{i6} (\vu - \vc{V})
\cdot \vc{n} |_{\Gamma_\tau} = 0 \ \mbox{for any}\ \tau \geq 0,
\eF
where $\vc{n}(t,x)$ denotes the unit outer normal vector to the boundary $\Gamma_t$.
Moreover, we study the problem with the full slip boundary conditions in
the form
\bFormula{b1} \left[ \tn{S} \vc{n} \right]\times \vc{n} = 0.
\eF
%$\vu$ and $\vc{V}$ denote the fluid and solid body velocities,
%respectively,
%and $\Gamma_\tau$ is the position of the boundary
%at a time $\tau$ with the outer normal vector $\vc{n}$.
%If
%$\zeta= 0$, we obtain the \emph{complete slip} while the
%asymptotic limit $\zeta \to \infty$ gives rise to the standard
%no-slip boundary conditions.
%Throughout this paper we work for clarity of presentation with a \emph{complete slip} condition, i.e. we set $\zeta = 0$.
% However it is easy to observe that the result of this paper hold also for $\zeta > 0$, see also Section \ref{d}.
%It is worth mentioning that $\zeta \to \infty$ gives rise to the standard no-slip boundary conditions

The heat flux satisfies the conservative boundary
conditions
\bFormula{b2} \vq\cdot \vc{n}=0 \mbox { for all }
t\in [0,T],\  x \in  {\Gamma_t}.
\eF

For the physical motivation of correct description of the fluid
boundary behavior, see Bul\' \i \v cek, M\' alek and Rajagopal
\cite{BMR1}, Priezjev and Troian \cite{PRTR} and the references
therein.

Finally we supplement our considered system with given initial data
	$$\varrho_0, \ (\varrho\vc{u})_0, \ \vt_0.$$

In comparison with \cite{KrMaNeWr}, we consider entropy \eqref{i1c} and energy \eqref{i1d} balance in place of thermal energy equation -- cf. \cite[(1.3)]{KrMaNeWr}. Even though strong solutions for both variants coincide as far as $\vt$ is bounded below away from zero, results concerning weak solutions cannot be simply transferred between them -- for example a relative entropy inequality which plays a significant role in weak-strong uniqueness and in various singular limits -- see \cite{BrKrMa,BumJa,FeNo6,FeNows}. Hence, it is reasonable to complete the theory by showing the existence of variational solution to system \eqref{i1a}--\eqref{i1d} provided $\Omega_t\subset \mathbb R^3$ is time-dependent. This is the main goal of this paper.

The existence theory for the barotropic Navier-Stokes system on \emph{fixed} spatial domains in the framework of weak solutions dates back to the seminal work by Lions \cite{LI4}, who worked with certain growth of pressure, and Feireisl et. al \cite{FNP} where the existence to a class of physically relevant pressure-density state equations was shown. These results were later extended to the full Navier-Stokes-Fourier system in \cite {EF70, EF71}, where the formulation with thermal energy equation is used.

Feireisl and Novotn\'y \cite{FeNo6} also proved the existence of weak solutions to the Navier-Stokes-Fourier system formulated with an entropy inequality instead of the thermal energy equation. This approach on the one hand allows for proving existence of solutions with more general hypotheses on the pressure on the other hand it requires the presence of a radiative term $a\vt^4$ in the pressure law.

The investigation of \emph{incompressible} fluids in {time dependent} domains started with a seminal paper of Ladyzhenskaya \cite{LAD2}, Fujita et al.\cite{FN}, see also \cite{Neust1,NeuPen1,NeuPen2,S} for more recent results in this direction.

\emph{Compressible} fluid flows in time dependent domains in barotropic case were examined in \cite{FeNeSt} for the \emph{no-slip} boundary conditions and in \cite{FKNNS} for the slip boundary conditions. As mentioned earlier, the existence of solution to the NSF system was the content of \cite{KrMaNeWr}. %We would like to emphasise that our approach is working for different restrictions of physical quantities.

In order to give a proof of existence of variational solutions, we proceed in the following way.
\begin{enumerate}

\item In order to be able to start our consideration of  the penalised system in a fixed time independent domain $B$ large enough, such that
\bFormula{set_B}
\overline{\Omega_t} \subset B \mbox{ for any }t\in [0,T]
\eF
we add to the momentum equation an extra term
 \bFormula{i8}
\frac{1}{\ep} \int_0^T \int_{\Gamma_t}  (\vu - \vc{V} ) \cdot \vc{n} \ \vph \cdot \vc{n} \, {\rm dS}_{x} \dt,\ \ep > 0 \ \mbox{small},
\eF
which was originally proposed by Stokes and Carey in \cite{StoCar}.
This penalizes the flux through the interface $\Gamma_t$ and allows to deal with the slip boundary conditions. Namely, as $\ep \to 0$, this additional term yields the boundary condition $(\vu - \vc{V})\cdot \vc{n} = 0$ on $\Gamma_t$. Consequently the large domain $(0,T) \times B$ becomes divided by an impermeable interface $\cup_{t\in (0,T)} \{t\} \times \Gamma_t$ to a \emph{fluid domain} $Q_T$ and a \emph{solid domain} $Q_T^c$. Next to handle the behaviour of the solution in the solid domain $Q_T^c$ we use a penalization scheme which is represented by parameters $\omega,\ \nu, \ \delta, \ \lambda$ and $\eta$.

%This extra term allows to consider the system in a large fixed time independent domain (denoted further by $B$) and penalizes the flux through the interface $\Gamma_t$. In the limit $\ep \rightarrow 0$ this term yields the boundary condition $(\vu - \vc{V})\cdot \vc{n} = 0$ on $\Gamma_t$, and thus the large domain $(0,T) \times B$ becomes divided by an impermeable interface $\cup_{t\in (0,T)} \{t\} \times \Gamma_t$ to a \emph{fluid domain} $Q_T$ and a \emph{solid domain} $Q_T^c$. In order to handle the behaviour of the solution in the solid domain $Q_T^c$ we use a penalization scheme which is represented by parameters $\omega,\ \nu, \ \delta$ and $\eta$.

\item In addition to  \eqref{i8}, we introduce \emph{variable} coefficients: the shear viscosity coefficient $\mu_\omega(t,x,\vt)$, the heat conductivity coefficient  $\kappa_{\nu}(t,x,\vt)$, and coefficient in a radiation counterpart of the pressure, the specific entropy, and internal energy function $a=a_\eta(t,x)$. All of them remain strictly positive in the fluid domain $Q_T$, but vanish in the solid domain $B \setminus \overline{Q_T^c}$ as $\omega,$ $\nu,$ and  $\eta$ converge to zero respectively.

%The coefficient $a:=a_\eta(x)$ which represents the radiative part (see \eqref{p1p}, \eqref{e1e} and \eqref{s1s}) is assumed to vanish on the solid domain. I.e. we assume that $a_\eta(t,x) = \left\{\begin{array}{l} a,\ a>0\ \mbox{ for } (t,x)\in Q_T\\ \eta a\ \mbox{ for } (t,x)\in Q_T^c.\end{array}\right.$ We use a letter $a$ for both coefficient and constant. Anyway, it is always clear which meaning is considered.
%
%\item In addition to (\ref{i8}), we introduce a \emph{variable}
%shear viscosity coefficient $\mu = \mu_\omega$, where $\mu_\omega$
%remains strictly positive in the fluid domain $Q_T$ but vanishes
%in the solid domain $Q_T^c$ as $\omega \to 0$.
%\item  We introduce the heat conductivity coefficient $\kappa_{\nu}(t,x,\vt)$ which
%remains strictly positive in the fluid domain $Q_T$ but vanishes
%in the solid domain $Q_T^c$ as $\nu \to 0$.
%For technical reasons this procedure needs to be done in two steps, more precisely we consider a step function $\kappa_\nu = \kappa$ in the fluid part and $\kappa_\nu = \nu\kappa$ in the solid part.
%This function is then mollified to get smooth $\kappa_{\nu,\xi}$, $\xi > 0$ with $\kappa_{\nu,\xi} \sil \kappa_\nu$ as $\xi \sil 0$.

\item We add a term $\lambda \vt^5$ into the energy balance and $\lambda \vt^4$ into the entropy balance. These terms yield a control over a temperature on the solid domain. The exact choice of power $\vt^5$ is not essential, however it is important that the power is larger than $\vt^4$.

\item Similarly to
the existence theory developed in \cite{FeNo6}, we introduce the
\emph{artificial pressure}
%, \emph{artificial energy}, \emph{artificial entropy}
related to coefficient $\delta$:
\[
p_\delta(\varrho,\vt ) = p(\varrho,\vt) + \delta \varrho^\beta ,\
\beta \geq 4,\ \delta > 0.
\]
This gives an extra information abut the density.
%\comment{Compare it pleas with Section 3.3.1 in \cite{FeNo6}. There is introduced also a penalization in $\mathbb S$, however, I think that this is not needed in our case. -- Venca}
\item Keeping $\ep,\,\eta,\,\omega,\,\nu,\,\lambda,$ and $\delta>0$ fixed, we solve the modified problem in a (bounded) reference domain $B\subset \mathbb{R}^3$ such that \eqref{set_B} is satisfied.
To this end, we adapt the existence theory for the compressible
Navier-Stokes-Fourier system with variable coefficients
developed in \cite{FeNo6}.

\item First, we take the initial density
$\varrho_0$ vanishing outside $\Omega_0$ and letting $\ep \to 0$
for fixed $\eta,\, \omega,\, \nu, \lambda,\, \delta > 0$ we obtain a ``two-fluid''
system where the density vanishes in the solid part $\left( (0,T)
\times B \right) \setminus Q_T$ of the reference domain.

%\item Passing with $\xi \sil 0$ we recover the system with jump in the heat conductivity coefficient and justify the choice of the boundary condition of the test function in the weak formulation of the thermal energy balance.
\item In order to get rid of the term on $Q_T^c$, we tend with remaining approximations to zero. We do it in the following sequence $\eta,\,\omega,\,\nu,\,\lambda,$ and $\delta$.

\end{enumerate}

The paper is organized as follows. In Section \ref{m} we present all assumptions and we complete system \eqref{i1a}--\eqref{i1d} by prescribing the boundary and initial conditions. We define a variational solution and we introduce the precise version of our main result.  Section \ref{p} is devoted to the penalization problem, we highlight all approximations and we discuss the existence of its solution. Finally, the proof of the main theorem is concluded in Section \ref{s} by performing appropriate limits.

\section{Preliminaries}
\label{m}

\subsection{Hypotheses}
\label{s:m1}

Motivated by \cite{FeNo6} we introduce the following set of assumptions:

The stress tensor $\tn{S}$ is determined by the standard Newton rheological law
\bFormula{i4}
\tn{S} (\vt, \Grad \vu) = \mu(\vt) \left( \Grad \vu +
\Grad^t \vu - \frac{2}{3} \Div \vu \tn{I} \right) + \eta (\vt) \Div \vu
\tn{I},\,\, \mu > 0,\ \eta \geq 0.
\eF
We assume the viscosity coefficients $\mu$ and $\eta$ are continuously differentiable functions of the absolute temperature, namely $\mu,\ \eta \in C^1[0,\infty)$ and satisfy
	\bFormula{mu_1}
	0< \underline{\mu} (1+ \vt) \leq \mu(\vt) \leq \overline{\mu} (1+\vt),
	\quad
	\sup_{\vt \in [0,\infty)} |\mu'(\vt)| \leq \overline{m},
	\eF
	\bFormula{eta_1}
	0 \leq \eta(\vt) \leq \overline{\eta} (1+\vt).
	\eF
%\comment{Later on I suppose that $\alpha =1$. The case $\alpha\leq1$ can be done but it requires some additional modifications, especially estimates of velocity. -- Venca}
The Fourier law for the heat flux $\vq$ has the following form:
 \bFormula {q1}
 \vq= -\kappa (\vt)\Grad\vt,
\eF
where the heat coefficient $\kappa$ can be decompose into two parts
	\bFormula{kappa_1}
	\kappa(\vt) = \kappa_M (\vt) + \kappa_R (\vt)
	\eF
where $\kappa_M, \ \kappa_R \in C^1[0,\infty)$ and
	\bFormula{kappa_2}
	0< \underline{\kappa_R} (1+ \vt^3) \leq \kappa_R (\vt) \leq \overline{\kappa_R}(1+ \vt^3),
	\eF
	\bFormula{kappa_3}
	0< \underline{\kappa_M}   (1+ \vt) \leq \kappa_M (\vt) \leq \overline{\kappa_M}  (1+ \vt).
	\eF
In the above formulas $\underline{\mu}$, $\overline{\mu}$, $\overline{m}$, $\overline{\eta}$, $\underline{\kappa_R}$, $\overline{\kappa_R}$, $\underline{\kappa_M} $, $\overline{\kappa_M}$ are positive  constants. Let us remark that the existence of solutions for the fixed domain can be obtain for more general  $\mu$ and $\kappa_M$ (see \cite{FeNo6}), namely with upper and lower growth described by function $(1+\vt^\alpha)$ with $\alpha \in(\frac{2}{5},1]$ instead of $(1+\vt)$. We believe our result can be extended to this more general one, however to simplify our consideration we assume $\alpha = 1$.

The entropy production rate $\sigma$ satisfies
	\bFormula{sigma1}
	\sigma \geq \frac{1}{\vt} \left( \tn{S} : \Grad \vu  -\frac{\vq}{\vt} \cdot \Grad \vt \right).
	\eF

The quantities $p$, $e$, and $s$ are continuously differentiable functions for positive values of $\vr$, $\vt$ and satisfy Gibbs' equation
	\bFormula{gibs}
	\vt D s(\vr,\vt) = D e(\vr,\vt) + p(\vr,\vt) D\left( \frac{1}{\vr}  \right) \mbox{ for all } \vr, \ \vt > 0.
	\eF

Further, we assume the following state equation for the pressure and the internal energy
\bFormula{p1p}
p(\vr,\vt) = p_M(\vr,\vt)+ p_R(\vt), \quad p_R(\vt) = \frac a3 \vt^4,\ a>0,
\eF
\bFormula{e1e}
e(\vr,\vt) = e_M(\vr,\vt) + e_R(\vr,\vt), \quad \vr e_R(\vr,\vt) = a\vt^4,
\eF
and
\bFormula{s1s}
s(\vr,\vt) = s_M(\vr,\vt) + s_R(\vr,\vt), \quad \vr s_R(\vr,\vt) = \frac 43 a\vt^3.
\eF
According to the hypothesis of thermodynamic  stability the molecular components satisfy
	\bFormula{pm_1}
	\frac{\partial p_M}{\partial \vr} >0 \mbox{ for all } \vr, \ \vt >0
	\eF
and
	\bFormula{em1}
	0<\frac{\partial e_M}{\partial \vt} \leq c \mbox{ for all }  \vr, \ \vt >0.
	\eF
Moreover
	\bFormula{em2}
	\lim_{\vt\to 0^+} e_M(\vr,\vt) = \underline{e}_M(\vr) > 0 \mbox{ for any fixed } \vr >0,
	\eF	
and
	\bFormula{em3}
	\left| \vr \frac{\partial e_M(\vr,\vt)}{\partial \vr} \right| \leq c e_M(\vr, \vt) \mbox{ for all } \vr, \ \vt >0.
	\eF
We suppose also that there is a function $P$ satisfying
	\bFormula{p2p}
	P \in C^1[0,\infty), \ P(0)=0,\  P'(0)>0,
	\eF
and two positive constants $0< \underline{Z} < \overline{Z}$ such that
	\bFormula{pm2}
	p_M(\vr,\vt) =  \vt^{\frac{5}{2}} P\left( \frac{\vr}{\vt^{\frac{3}{2}}} \right)
	\mbox{ whenever } 0< \vr \leq \underline{Z} \vt^{\frac{3}{2}},\mbox{ or, } \vr > \overline{Z} \vt^{\frac{3}{2}}
	\eF
and
	\bFormula{pm3}
	p_M(\vr,\vt) = \frac{2}{3} \vr e_M (\vr,\vt) \mbox{ for } \vr > \overline{Z}\vt^{\frac{3}{2}}.
	\eF
%\comment{The assumptions \eqref{pm_1} -- \eqref{pm3} should be checked. Do they match with \cite{FeNo6}? -- Venca}

Finally, the problem \eqref{i1a}--\eqref{i1d} is supplemented by the
initial conditions
\bFormula{i7_1}
\vr(0, \cdot) = \vr_0 \in L^{\frac{5}{3}}(\Omega_0),\quad \vr_0 \geq 0, \quad \vr_0 \not\equiv 0,\quad \vr_0|_{\mathbb{R}^3 \setminus \Omega_0} = 0,
\eF
\bFormula{i7_2}
(\vr\vu) (0, \cdot) = (\vr \vu)_0, \ (\vr \vu)_0 = 0 \mbox{ a.e.  on the set }
\{ \Omega_0 \ | \ \vr_0(x) = 0 \},\quad
\int_{\Omega_0} \frac{| (\vr\vu)_0 |^2}{\vr_0}\, \dx < \infty ,
\eF
\bFormula{i7_3}
\vt_0 >0 \mbox{ a.e. in }\Omega_0, %{\color{red} \vt_0 \in L^\infty(\Omega_0), \ \vt_0 \leq \underline \vt >0 \mbox{ in }\Omega_0, }
 \qquad   (\vr s)_0 = \vr_0 s(\vr_0, \vt_0) \in L^1 (\Omega_0),
\eF
\bFormula{i7_4}
E_0 = \int_{\Omega_0} \left( \frac{1}{2\vr_0} |(\vr \vu)_0 |^2  + \vr_0 e(\vr_0,\vt_0)  \right)dx < \infty .
\eF
%\comment{red part; this remains probably from the previous paper, it is not needed in \cite{FeNo6} -- Aneta}
%%%%%%
%%%%%%
\subsection{Weak formulation, main result}

%In the weak formulation, it is convenient to consider the continuity equation \eqref{i1a} in the whole physical space $\mathbb{R}^3$ provided the density
%$\vr$ is extended to be equal to zero outside the fluid domain, specifically
%\bFormula{m1}
%\int_{\Omega_\tau} \vr \varphi (\tau, \cdot) \ \dx - \int_{\Omega_0} \vr_0 \varphi (0, \cdot) \ \dx =
%\int_0^\tau \int_{ \Omega_t} \left( \vr \partial_t \varphi + \vr \vu \cdot \Grad \varphi \right) \ \dxdt
%\eF
%for any $\tau \in [0,T]$ and any test function $\varphi \in \DC([0,T] \times \mathbb{R}^3)$.
In the weak formulation, equation (\ref{i1a}) is supposed to be fulfilled in the sense of
renormalized solutions introduced by DiPerna and Lions \cite{DL}:
\bFormula{m2}
\int_0^T \int_{\Omega_t} \vr B(\vr) ( \partial_t \varphi + \vu \cdot \Grad \varphi )\, \dxdt =
\int_0^T \int_{\Omega_t} b(\vr)  \Div \vu \varphi\, \dxdt - \int_{\Omega_0}  \vr_0 B(\vr_0)  \varphi (0)\, \dx
\eF
 for any $\varphi \in C^1_c([0,T) \times \mathbb{R}^3)$, and any  $b \in L^\infty \cap C [0, \infty)$ such that { $b(0)=0$} and
 $B(\vr) = B(1) + \int_1^\vr \frac{b(z)}{z^2} {\rm d}z.$
Of course, we suppose that $\vr \geq 0$ a.e. in $(0,T) \times \mathbb{R}^3$.
% \comment{A comment to one of previous doubts: It seems to me that we need $b(0)=0$ to kill the first term on the right hand side outside of fluid domain, so I have added this restriction back. - Aneta}

The momentum equation \eqref{i1b} is replaced by the following integral identity
\bFormula{m3}
\int_0^T \int_{\Omega_t} \left( \vr \vu \cdot \partial_t \vph + \vr [\vu \otimes \vu] : \Grad \vph + p(\vr,\vt) \Div \vph
- \tn{S} (\vt, \Grad \vu) : \Grad \vph\right) \dxdt
= - \int_{\Omega_0} (\vr \vu)_0 \cdot \vph (0, \cdot) \ \dx,
 \eF
which should be fulfilled for any test function $\vph \in C^1_c (\overline{Q_T} ; \mathbb{R}^3)$
such that $\varphi(T,\cdot) =0$ and
\bFormula{m4}
\vph \cdot \vc{n}|_{\Gamma_\tau} = 0 \ \mbox{for any} \ \tau \in [0,T].
\eF

The impermeability condition (\ref{i6}) is satisfied in the sense of traces, specifically,
\bFormula{m50}
\vu,\nabla_x\vu \in L^2(Q_T; \mathbb{R}^3) \ \mbox{and}\ (\vu - \vc{V}) \cdot \vc{n}  (\tau , \cdot)|_{\Gamma_\tau}  = 0 \ \mbox{for a.a.}\ \tau \in [0,T].
\eF

The entropy inequality
	\begin{multline}
	\label{m5}
	\int_0^T \int_{\Omega_t} \vr s (\partial_t \varphi + \vu \cdot \Grad \varphi ) \dxdt
	- \int_0^T \int_{\Omega_t} \frac{\kappa(\vt) \Grad \vt \cdot \Grad \varphi }{\vt} \dxdt
	+ \int_0^T \int_{\Omega_t} \frac{\varphi}{\vt} \left( \tn{S} : \Grad \vu + \frac{\kappa(\vt) | \Grad \vt |^2}{\vt} \right)\\
	\leq - \int_{\Omega_0} (\vr s )_0 \varphi (0) \dx
	\end{multline}
holds for all $\varphi \in C^1_c (\overline{Q_T})$ such that $\varphi(T,\cdot) =0$ and $\varphi \geq 0$.

Finally, we assume the following energy inequality
\begin{multline}\label{m6}
\int_{\Omega_\tau}\left(\frac12 \vr |\vu|^2 + \vr e\right)(\tau,\cdot) \dx\leq \int_{\Omega_0} \left(\frac12 \frac{(\vr \vu)_0^2}{\vr_0} + \vr_0 e_0 - (\vr \vu)_0\cdot \vc{V}(0)\right)\dx\\
-\int_0^\tau\int_{\Omega_t}\left( \vr(\vu\otimes\vu):\Grad \vc{V} + p\, \Div \vc{V}  -\mathbb S:\Grad \vc{V} + \vr \vu\cdot \partial_t \vc{V}\right) \dxdt  + \int_{\Omega_t}\vr \vu \cdot \vc{V}(\tau,\cdot)\dx
\end{multline}
holds for a.a. $\tau\in(0,T)$.
%Balance of total energy
	%\bFormula{m6}
	%\int_0^T  \int_{\Omega_t} \left( \frac{1}{2} \vr |\vu|^2 + \vr e    \right) \partial_t \psi  \dx\dt =
	%-  \int_{\Omega_0}  \left( \frac{1}{2} \frac{ ( \vr \vu)_0^2}{\vr_0}  + \vr_0 e_0    \right)  \psi (0) \dx
	%+ \int_0^T  \int_{\Omega_t} \vr \vu \cdot \partial_t ( \vc{V}  \psi) +  \vr \vu \otimes \vu : \Grad \vc{V} \psi
	%\eF
	%\[
	%+ p(\vr) \Div \vc{V} \psi  - \tn{S} : \Grad \vc{V}\psi  - \vr\vf \cdot (\vu -\vc{V})\psi \dx\dt + \int_{\Omega_0} (\vr\vu)_0 \cdot \vc{V}(0) \psi (0) \dx
	%\]
%is fulfilled for all $\psi \in C^1_c ([0,T))$.

%\comment{Compare it with \eqref{final.ene} (lack of $\psi$). Equality or inequality? -- Venca}

\begin{Definition}\label{d:WS}
We say that the trio $(\vr,\vu,\vt)$ is a \emph{variational solution} of problem \eqref{i1a}--\eqref{i1d} with boundary conditions \eqref{i6}--\eqref{b2} and initial conditions \eqref{i7_1}--\eqref{i7_4} if
\begin{itemize}
\item $\vr \in L^\infty(0,T; L^{\frac 53}(\mathbb R^3))$, $\vr \geq 0$,
	$\vr \in L^q(Q_T)$ for certain $q>\frac{5}{3}$,
\item $\vu,\, \nabla \vu \in L^2(Q_T)$, $\vr \vu \in L^\infty(0,T;L^1(\mathbb R^3))$,
\item  $\vt>0$ a.a. on $Q_T$,
$\vt \in L^\infty((0,T);L^4 (\R^3))$, $\vt,\,\nabla \vt \in L^2(Q_T),$ and  $\log\vt,\,\nabla \log\vt \in L^2(Q_T),$
\item $\vr s,\ \vr s \vu, \ \frac{\vq}{\vt} \in L^1(Q_T)$,
\item relations \eqref{m2}--\eqref{m6} are satisfied.
\end{itemize}
\end{Definition}

\begin{Remark}
In contrast to \cite{FeNo6} we consider energy inequality rather than energy equation. Although it seems that we are losing a lot of information our definition of weak solution is still sufficient. Namely, if the above defined weak solution is smooth enough it will be a strong one. For a justification of this fact we refer to \cite[Section 1.2]{poul} .
\end{Remark}

At this stage, we are ready to state the main result of the  present paper:

\bTheorem{m1}
Let $\Omega_0 \subset \mathbb{R}^3$ be a bounded domain of class $C^{2 + \nu}$ with some $\nu >0$, and let $\vc{V} \in C^1([0,T]; C^{3}_c (\mathbb{R}^3;\mathbb{R}^3))$ be given and satisfy \eqref{V_0}, \eqref{Vselenoidal}. Assume that hypothesis \eqref{i4}--\eqref{pm3} are satisfied.
%
%Furthermore let the initial data fulfill
%\[
%\vr_0 \in L^\gamma (\mathbb{R}^3),\ \vr_0 \geq 0, \ \vr_0 \not\equiv 0,\ \vr_0|_{\mathbb{R}^3 \setminus \Omega_0} = 0,\
%(\vr \vu)_0 = 0 \ \mbox{a.a. on the set} \ \{ \vr_0 = 0 \} ,\ \int_{\Omega_0} \frac{1}{\vr_0} |(\vr \vu)_0 |^2 \ \dx < \infty
%\]
%and
%\[
%\vt_0  \in L^\infty(\Omega_0),\quad \vt_0 \geq \underline{\vt}  > 0 \quad \mbox{ on } \Omega_0.
%\]

Then the problem \eqref{i1a}--\eqref{i1d} with boundary conditions \eqref{b1}, \eqref{b2} and initial conditions \eqref{i7_1}--\eqref{i7_4} admits a variational solution in the sense of Definition \ref{d:WS} on any finite time interval $(0,T)$.
\eT
%\comment{Do we really need such smooth domain? -- Venca}
The rest of the paper is devoted to the proof of Theorem \ref{Tm1}.

\section{Approximate problem}
\label{p}

\subsection{Penalized problem - weak formulation}

Choosing $R > 0$ such that
$$
\vc{V} |_{[0,T] \times \{ |x| > R \} } = 0 ,\ \ \Ov{ \Omega }_0 \subset \{ |x| < R \}
$$
we take the reference domain $B = \{ |x| < 2R \}$.

The shear viscosity coefficients  $\mu_\omega$ and $\eta_\omega$ are taken such that
\bFormula{p2}
\mu_\omega(\vt,\cdot) \in \DC ([0,T] \times \mathbb{R}^3),\ 0 < \omega\mu(\vt) \leq \mu_\omega (\vt,t,x) \leq \mu(\vt) \ \mbox{in}\ [0,T] \times B, \
\mu_{\omega}(\vt,\tau, \cdot)|_{ \Omega_\tau } = \mu(\vt) \ \mbox{for any} \ \ \tau \in [0,T]
\eF
\bFormula{p2_eta}
\eta_\omega(\vt,\cdot) \in \DC ([0,T] \times \mathbb{R}^3),\ 0< \omega \eta(\vt) \leq \eta_\omega (\vt,t,x) \leq \eta(\vt) \ \mbox{in}\ [0,T] \times B, \
\eta_{\omega}(\vt,\tau, \cdot)|_{ \Omega_\tau } = \eta(\vt) \ \mbox{for any} \ \ \tau \in [0,T]
\eF
and
\bFormula{p2_2}
 \mu_\omega, \eta_\omega \to 0 \mbox{ a.e. in } ((0,T) \times B )\setminus Q_T \qquad \mbox{ as } \omega \sil 0.
\eF

We introduce a variable heat conductivity coefficient as follows:
\bFormula{k_nu}
\kappa_\nu(\vt,t,x) = \chi_{\nu}(t,x) \kappa(\vt), \mbox{ where }
\chi_\nu = 1  \mbox{ in } Q_T  \mbox{ and } \chi_\nu = \nu  \mbox{ in } ((0,T) \times B) \setminus Q_T.
\eF

Similarly we introduce a variable coefficient $a:=a_\eta(t,x)$ which represents the radiative part of pressure, internal energy and entropy (see \eqref{p1p}, \eqref{e1e}, and \eqref{s1s}). Namely, we assume that
\bFormula{a_delta}
a_\eta(t,x) = \chi_{\eta}(t,x)a, \mbox{ where } a>0 \mbox{ and } \chi_\eta = 1  \mbox{ in } Q_T  \mbox{ and } \chi_\eta = \eta  \mbox{ in } ((0,T) \times B \setminus Q_T.
\eF
We use a letter $a$ for both coefficient and constant. Anyway, it is always clear which meaning is considered.
%\comment{I erased the mollification. Hope it does not bring any problem later on. We have to take care.

%On the other hand it seems that $\kappa_\nu\to 0$ a.e. in $B\setminus \Omega_t$ is sufficient so we can have continuous $\kappa_\nu$ -- Venca}
Moreover being motivated by the approximation for the existence theory we follow  \cite{FeNo6} and set:
\bFormula{p_delta}
p_{\eta,\delta}(\varrho,\vt ) = p_M(\varrho,\vt) + \frac{a_\eta}{3}\vt^4 + \delta \varrho^\beta ,\
\beta \geq 4,\ \delta > 0,
\eF
\bFormula{e_s_delta}
e_\eta(\varrho, \vt) = e_{M}(\varrho,\vt) + a_\eta\frac{\vt^4}{\vr},
\quad\quad
s_\eta(\varrho,\vt) =  s_M(\vr,\vt) +  \frac 43 a_\eta \frac{\vt^3}{\varrho}.
\eF
%\bFormula{k_delta}
%\kappa_{\nu} =  \chi_{\nu}(t,x) \kappa(\vr,\vt) + \frac{\delta}{\vt}).
%\eF

Finally, let $\vr_0$, $(\vr\vu)_0$ and $\vt_0$ be initial conditions as specified in Theorem \ref{Tm1}. We define modified initial data $\vr_{0,\delta}$, $(\vr\vu)_{0,\delta}$ and $\vt_{0,\delta}$ so that %for $\delta > 0$ we define modified initial data so that
\bFormula{data1}
\vr_{0, \delta} \geq 0, \ \vr_{0,\delta} \not\equiv 0,\ \vr_{0, \delta}|_{\mathbb{R}^3 \setminus \Omega_0} = 0,\ \int_{B}
\left( \vr_{0, \delta}^{\frac{5}{3}} + \delta \vr_{0, \delta}^\beta \right) \dx \leq c,\ \vr_{0,\delta}\sil \vr_0 \ \mbox{in } L^{\frac{5}{3}}(B), \ |\{\vr_{0,\delta}<\vr_0\}|\sil 0,
\eF
\bFormula{data2}
%(\vr \vu)_0 = (\vr \vu)_{0, \delta},\ (\vr \vu)_{0, \delta} = 0 \ \mbox{a.a. on the set} \ \{ \vr_{0, \delta} = 0 \} ,\ \int_{\Omega_0} \frac{1}{\vr_{0, \delta}} |(\vr \vu)_{0,\delta} |^2 \ \dx \leq c
(\vr\vu)_{0,\delta} = \left\{
\begin{array}{ll}
(\vr\vu)_{0}\ & \mbox{if}\ \vr_{0,\delta}\geq \vr_0,\\
0 & \mbox{otherwise} ,
\end{array}
\right.
\eF
	\bFormula{data3}
	\quad 0< \underline{\vt} \leq \vt_{0,\delta}
    \mbox{ and }\vt_{0,\delta} \in L^\infty(B) \cap C^{2+\nu} (B).
	\eF
Moreover
$$
\int_{\Omega_0} \vr_{0,\delta}e(\vr_{0,\delta},\vt_{0,\delta}) \dx\to \int_{\Omega_0} \vr_0 e(\vr_0,\vt_0)\dx
$$
and
$$
\vr_{0,\delta}s(\vr_{0,\delta},\vt_{0,\delta})\to \vr_0s(\vr_0,\vt_0) \ \mbox{weakly in}\ L^1(\Omega_0).
$$

%\comment{The last convergence is maybe not sufficient. Again, the convergence of entropy should be somehow engaged. This should follow from the very last section. -- Venca}

Now we are ready to state the weak formulation of the \emph{penalized problem}.

 Again, we consider $\vr,\vu$ to be zero outside of $(0,T) \times B$. The weak (renormalized) formulation of the continuity equation reads as
\bFormula{p3}
\int_0^T \int_B \vr B(\vr) ( \partial_t \varphi + \vu \cdot \Grad \varphi ) \, \dxdt =
\int_0^T \int_B b(\vr)  \Div \vu \varphi \dxdt - \int_B  \vr_{0,\delta} B(\vr_{0,\delta})  \varphi (0)\, \dx
\eF
 for any $\varphi \in C^1_c([0,T) \times \mathbb{R}^3)$, and any  $b \in L^\infty \cap C [0, \infty)$ such that $b(0)=0$ and
 $B(\vr) = B(1) + \int_1^\vr \frac{b(z)}{z^2} {\rm d}z.$
 The momentum equation is represented by the family of integral identities
\bFormula{p4}
 \int_0^T \int_{B} \left( \vr \vu \cdot \partial_t \vph + \vr [\vu \otimes \vu] : \Grad \vph + p_{\eta,\delta}(\vr,\vt) \Div \vph
- \tn{S}_\omega : \Grad \vph \right) \, \dxdt
 - \frac{1}{\ep} \int_0^\tau \int_{ \Gamma_t } \left( ( \vu - \vc{V} ) \cdot \vc{n} \ \vph \cdot \vc{n} \right) \, {\rm dS}_{x} \dt
\eF
\[
=  -  \int_{B} (\vr \vu)_{0,\delta} \cdot \vph (0, \cdot) \, \dx
\]
\[
\mbox{ with }\quad \tn{S}_\omega = {\mu_\omega}(\vt,t,x)
 \left( \Grad \vu + \Grad^t \vu - \frac{2}{3} \Div \vu \tn{I} \right) + \eta_\omega(\vt,t,x) \Div \vu \mathbb I
\]
for any test function $\vph \in \DC([0,T) \times B ; \mathbb{R}^3)$.
	\bFormula{p5}
	 \int_0^T \int_{B}   \vr s_\eta (\partial_t \varphi + \vu \cdot \Grad \varphi )\, \dxdt
	- \int_0^T \int_{B}\frac{\kappa_{\nu}(\vt,t,x) \Grad \vt \cdot \Grad \varphi }{\vt}\, \dxdt
		\eF
		\[
	+  \int_0^T \int_{B} \frac{\varphi}{\vt} \left( \tn{S}_\omega : \Grad \vu
	+ \frac{\kappa_{\nu}(\vt,t,x) | \Grad \vt |^2}{\vt}   \right)\, \dxdt - \int_0^T \int_{B} \lambda \vt^4 \dxdt
	\leq - \int_B (\vr s )_{0,\delta,\eta} \varphi (0) \dx
	\]
for all $\varphi \in C^1_c ([0,T) \times \overline{B})$, $\varphi \geq 0$, where $(\vr s)_{0,\delta,\eta} := \vr_{0,\delta}s_\eta(\vr_{0,\delta},\vt_{0,\delta})$ and
	$$\sigma_{\omega,\nu} \geq  \frac{1}{\vt}\left(\mathbb S_\omega : \Grad \vu + \frac{\kappa_{\nu} (\vt,t,x) |\Grad \vt|^2}{\vt} \right). $$
%\comment{What is $\xi$? -- Venca}

The solution to penalized problem should satisfy the energy equation

	\begin{multline}\label{p6}
	\int_0^T  \int_B \left( \frac{1}{2} \vr |\vu|^2 + \vr e_\eta  + \frac{\delta}{\beta -1} \vr^\beta  \right) \partial_t \psi -\lambda \vt^5 \psi\,  \dxdt =
	\frac{1}{\ep} \int_0^T \int_{ \Gamma_t } ( \vu - \vc{V} ) \cdot \vc{n} \ \vu \cdot \vc{n} \psi\, {\rm dS}_{x} \, \dt \\
	-   \int_B  \left( \frac{1}{2} \frac{ ( \vr \vu)_{0,\delta}^2}{\vr_{0,\delta}}  + \vr_{0,\delta} e_{0,\delta, \eta}
	+ \frac{\delta}{\beta -1} \vr_{0,\delta}^\beta
	  \right)  \psi (0) \dx
	\end{multline}
	for all $\psi \in C^1_c ([0,T))$. Here we denoted $e_{0,\delta,\eta} := e_\eta(\vr_{0,\delta},\vt_{0,\delta})$. However, we will rather work with the following modification which can be obtained from \eqref{p4} and \eqref{p6}:
	\bFormula{p7}
	\begin{split}
	\int_0^T  \int_B & \left( \frac{1}{2} \vr |\vu|^2 + \vr e_\eta + \frac{\delta}{\beta -1} \vr^\beta    \right) \partial_t \psi   - \lambda \vt^5\psi\dxdt
	- \frac{1}{\ep} \int_0^T \int_{ \Gamma_t }  \left| ( \vu - \vc{V} ) \cdot \vc{n} \right|^2\psi\, {\rm dS}_{x}  \dt
	\\ &
	= -  \int_B  \left( \frac{1}{2} \frac{ ( \vr \vu)_{0,\delta}^2}{\vr_{0,\delta} } + \vr_{0,\delta} e_{0,\delta,\eta}  + \frac{\delta}{\beta -1} \vr^\beta_{0,\delta}   -  (\vr \vu )_{0,\delta} \cdot \vc{V}(0)  \right)  \psi (0) \dx
	\\ &
	-  \int_0^T \int_B (
	  \tn{S}_\omega : \Grad \vc{V} \psi
	 - \vr \vu \cdot \partial_t (\vc{V} \psi) - \vr \vu \otimes \vu : \Grad \vc{V} \psi -
p_{\eta,\delta} \Div \vc{V} \psi)  \dxdt
	\end{split}
	\eF
for all $\psi \in C^1_c ([0,T))$.
%\comment{Do we get for the penalized problem the energy equality or an energy inequality? -- Venca}

\begin{Definition}\label{d:WSpen}
Let $\ep,\,\eta,\, \omega,\,\nu,\,\lambda$, and $\delta$ be positive parameters and let $\beta > 4$.
We say that a trio $(\vr,\vu,\vt)$ is a variational solution to the penalized problem with initial data \eqref{data1}--\eqref{data3} if
\begin{itemize}
\item $\varrho \in L^\infty((0,T); L^{\frac{5}{3}}(\R^3)) \cap  L^\infty((0,T); L^\beta(\R^3))$, $\vr \geq 0$,
	$\vr \in L^q((0,T)\times B)$ for certain $q>\beta$,
\item $\vu,\, \nabla \vu \in L^2((0,T) \times B)$, $\vr \vu \in L^\infty(0,T;L^1(\mathbb R^3))$,
\item  $\vt>0$ a.a. on $(0,T)\times B$,
$\vt \in L^\infty((0,T);L^4 (B))$, $\vt,\,\nabla \vt \in L^2((0,T) \times B),$ and  $\log\vt,\,\nabla \log\vt \in L^2((0,T) \times B),$
\item $\vr s,\ \vr s \vu, \ \frac{\vq}{\vt} \in L^1((0,T)\times B)$,
\item relations \eqref{p3}--\eqref{p5}, \eqref{p7} are satisfied.

 \end{itemize}

\end{Definition}

The choice of the no-slip boundary condition $\vu|_{\partial B} = 0$ is not essential here as the proof follows the one presented in \cite{FeNo6}. We have the following existence theorem concerning weak solutions to the penalized problem.

\bTheorem{m2}
Let $\vc{V} \in C^1([0,T]; C^{3}_c (\mathbb{R}^3;\mathbb{R}^3))$ be given. Let $\beta > 4$. Let thermodynamical functions and coefficients satisfy \eqref{p2}--\eqref{e_s_delta} with assumptions on $p$, $e$, $s$, $\mu$, $\eta$, $\kappa$ as in Theorem~\ref{Tm1}. Let initial data satisfy \eqref{data1}--\eqref{data3}.
Finally, let  $\ep,\, \eta,\, \omega,\, \nu,\, \lambda,\, \delta > 0$.

Then the penalized problem admits a variational  solution on any time interval $(0,T)$ in the sense specified by Definition \ref{d:WSpen}.

\eT
%\comment{The proof below has to be changed}
\begin{proof}
As mentioned above, the proof itself does not differ from \cite[Section 3]{FeNo6} so we present just a short description. It is necessary to regularize the continuity equation by a viscous term $\Delta \varrho$ and to add appropriate terms to the momentum and energy equation. The starting point of considerations at the level of Galerkin approximations is then the internal energy equation \cite[Equation (3.55)]{FeNo6} instead of the entropy balance. However we need to accommodate two additional difficulties, namely
\begin{enumerate}
\item The term $\frac{1}{\ep} \int_0^\tau \int_{ \Gamma_t } \left( ( \vu - \vc{V} ) \cdot \vc{n} \ \vph \cdot \vc{n} \right) \, {\rm dS}_{x} \dt$ in \eqref{p4} and $\frac 1\varepsilon \int_0^T\int_{\Gamma_t} ({\bf u} - {\bf V})\cdot {\bf n}\ {\bf u} \cdot {\bf n} \psi {\rm d}S_x\ {\rm d}t$ in \eqref{p6}. These terms do not cause much trouble as each can be treated as a "compact" perturbation.
\item The jumps in functions $\kappa_\nu(\vt,t,x)$ and $a_\eta(t,x)$. To this end we first introduce mollifications of the jump function $\chi_A(t,x)$ where stands for $\nu$ or $\eta$ (see \eqref{k_nu}, \eqref{a_delta}). By $\chi_A^\alpha(t,x)$ we denote smooth function with value $1$ on $Q_T$ and with value $A$ on $B \setminus Q_T^\alpha$, where $Q_T^\alpha$ denotes the $\alpha$-neighborhood of $Q_T$ in spacetime. The terms related to the parameter $\eta$ are treated in the straightforward way. The terms related to the parameter $\nu$ requires a little more attention.

Since we want the leading term in the internal energy equation \cite[(3.55)]{FeNo6} to be Laplacian, we need to add to this equation a term of the form $\Div(\K(\vt)\nabla_x\chi_\nu^\alpha)$, where $\K(\vt) = \int_1^\vt \kappa(z) dz$. This way we obtain
\begin{align}\label{eq:internal}
&\partial_t(\vr e_\eta(\vr,\vt)) + \Div(\vr e_\eta(\vr,\vt)\vu) - \Div \nabla_x(\chi_\nu^\alpha(t,x)\K(\vt))  \\ \nonumber
& \qquad = \tn{S}_\omega : \nabla_x\vu - p_{\eta,\delta}(\vr,\vt)\Div\vu - \Div(\K(\vt)\nabla_x\chi_\nu^\alpha(t,x)) - \lambda\vt^5.
\end{align}

We can then follow the theory in \cite[Section 3.4.2]{FeNo6} to deduce existence of strong solutions to equation \eqref{eq:internal} and when passing to the limit in the Galerkin approximations switch to the entropy balance and global total energy balance. Then we pass to the limit with the artificial viscosity parameter as in \cite[Section 3]{FeNo6}.

As a final step in the proof of Theorem \ref{Tm2} we need to pass with $\alpha$ to zero. Note that
\begin{equation}\label{eq:alpha}
\chi_A^\alpha \rightarrow \chi_A \qquad \text{strongly in } L^p(B) \text{ for any } p < \infty.
\end{equation}
The limit passage with $\alpha$ is straightforward in the terms related to parameter $\eta$, in particular due to the presence of term $\lambda \vt^5$ providing uniform bounds of high power of the temperature on $B$.

Concerning the terms related to $\nu$, there are two of them in the entropy balance \eqref{p5}, which we need to take care of. First, there is the nonnegative term $\int_0^T \int_B \chi_\nu^\alpha \kappa(\vt)\abs{\nabla_x\vt}^2\vt^{-2}\varphi \dxdt$. Due to the nonnegativity of the functions appearing in this term we can use the inequality $\chi_\nu \leq \chi_\nu^\alpha$ and weak lower semicontinuity to pass to the limit in the same way as we later explain in \eqref{s10}.

Since we have at hand the same apriori estimates as we work with in the first limit passage in the proof of Theorem \ref{Tm1}, we know in particular \eqref{p15_10} which together with \eqref{eq:alpha} is enough to pass to the limit with $\alpha$ in the term $\int_0^T \int_B \chi_\nu^\alpha \kappa(\vt)\vt^{-1}\nabla_x\vt\cdot\nabla_x\varphi \dxdt$.

%This can be solved at the beginning by solving the appropriate quasilinear parabolic system - see . It is sufficient to search for a solution $\vartheta$ fulfilling additionally $\nabla\vartheta\cdot \vc{n}|_{\Gamma_t} = 0$. The comparison principle as well as its corollaries work in the same way and the smooth solutions constructed in means of \cite[Theorem 10.24]{FeNo6} can should be treated separately on $\Omega_t$ and $B\setminus\Omega_t$.
\end{enumerate}
%Since there appears no other serious difficulty, the viscous approximation may vanish similarly as in \cite{FeNo6} and we immediately get the existence of the demanded solution.
\end{proof}

%In addition, since $\beta > 4$, the density is square integrable and we may use the regularization technique of DiPerna and Lions \cite{DL}
%to deduce the renormalized version of (\ref{p3}), namely
%\bFormula{p6a}
%\int_{B} b(\vr) \varphi (\tau, \cdot) \, \dx - \int_{B} b(\vr_{0,\delta}) \varphi (0, \cdot) \, \dx =
%\int_0^\tau \int_{B} \left( b(\vr) \partial_t \varphi + b(\vr) \vu \cdot \Grad \varphi +
%\left( b(\vr)  - b'(\vr) \vr \right) \Div \vu \varphi \right) \,\dxdt
%\eF
%for any $\varphi$ and $b$ as in (\ref{m2}).

%\comment{Possible problem coming from non-smooth $\kappa$? Need check. -- Venca}
%\end{proof}

\subsection{Modified energy inequality and uniform bounds}\label{s:bounds}

We use $\varphi(t,x) \equiv 1 \cdot \psi_\xi(t)$ where $\psi_\xi\in C^1_c([0,T))$ is non-increasing function which fulfills
\begin{equation}\label{psi_zeta}
\psi_\xi(t) = \left\{\begin{array}{l}1\ \mbox{for } t<\tau-\xi\\ 0\ \mbox{for }t\geq \tau\end{array}\right.\quad \mbox{for some }\tau\in(0,T)\ \mbox{and arbitrary }\xi>0
\end{equation}
 as a test function in \eqref{p5} in order to derive (after passing with $\xi \to 0$)
$$
-\int_{B} \vr s_\eta(\vr,\vt)(\tau,\cdot)\dx + \int_0^\tau\int_B \frac{1}{\vt}\left(\mathbb S_\omega : \Grad \vu + \frac{\kappa_{\nu} (\vt) |\Grad \vt|^2}{\vt} \right) \dxdt - \int_0^\tau \int_B \lambda \vt^4 \dxdt \leq -\int_B (\vr s)_{0,\delta,\eta} \dx.
$$
 The same test function applied in \eqref{p7} implies
\begin{multline*}
\int_B \left(\frac 12 \vr |\vu|^2 + \vr e_\eta   +  \frac{\delta}{\beta -1} \vr^\beta  \right)(\tau,\cdot) \dx + \frac 1\varepsilon \int_0^\tau\int_{\Gamma_t} |(\vu-\vc{V})\cdot \vc{n}|^2 \,{\rm d}S_x{\rm d}t
+ \int_0^\tau \int_B \lambda \vt^5 \dxdt
\\
=  \int_0^\tau\int_B\left( \mathbb S_{\omega} : \Grad \vc{V}- \vr \vu \cdot \vc{V}_t - \vr (\vu\otimes \vu):\Grad \vc{V} - p_{\eta,\delta}(\vr, \vt) \Div \vc{V} \right)\dxdt\\
 +\int_B \vr \vu \cdot \vc{V}(\tau,\cdot) \dx + \int_B \left(\frac 12 \frac{(\vr\vu)_{0,\delta}^2}{\vr_{0,\delta}} + \vr_{0,\delta} e_{0,\delta,\eta} +  \frac{\delta}{\beta -1} \vr^\beta_{0,\delta}   - (\vr \vu)_{0,\delta} \cdot  \vc{V}(0)\right) \dx
\end{multline*}
for almost all $\tau \in (0,T)$ and we deduce that
\begin{equation}\label{ae1}
\begin{split}
&\int_{B}\left(\frac 12 \vr |\vu|^2 + H_{1,\eta}(\vr,\vt) +   \frac{\delta}{\beta -1} \vr^\beta  \right)(\tau,\cdot)\ {\rm d}x
\\ &  + \frac 1\varepsilon \int_0^\tau \int_{\Gamma_t} |(\vu - \vc{V})\cdot \vc{n}|^2  \,{\rm d}S_x{\rm d}t+ \int_{0}^\tau \int_B \frac 1\vt \left(\mathbb S_\omega :\Grad \vu + \frac{\kappa_{\nu}(\vt)|\Grad \vt|^2}\vt \right) \dxdt
+ \int_0^\tau \int_B  \lambda\vt^5 \dxdt
\\
& \leq \int_0^\tau\int_B\left( \mathbb S_{\omega} : \Grad \vc{V} - \vr (\vu\otimes \vu):\Grad \vc{V} - \vr \vu \cdot \vc{V}_t - p_{\eta,\delta}(\vr, \vt) \Div \vc{V}  + \lambda \vt^4\right) \dxdt
\\
& + \int_B \vr \vu \cdot \vc{V}(\tau,\cdot)\dx
 + \int_B \left(\frac 12 \frac{(\vr \vu)_{0,\delta}^2}{\vr_{0,\delta}}  + H_{1,\eta}(\vr_{0,\delta},\vt_{0,\delta}) +  \frac{\delta}{\beta -1} \vr_{0,\delta}^\beta  - (\vr \vu)_{0,\delta} \cdot \vc{V}(0,\cdot)\right) \dx,
\end{split}
\end{equation}
%\comment{Please, check all the signs on the right hand side since I am not sure. I know they do not matter, however it would be fantastic to have them correct. -- Venca}
for almost all $\tau \in (0,T)$, where
$$H_{1,\eta}(\vr,\vt) = \varrho(e_\eta(\varrho,\vt) - s_\eta (\varrho,\vt))$$
(see \cite[Chapter 2.2.3]{FeNo6})
is a Helmholtz function.
%%and its approximations are defined as
%%$$
%%H_{1,\delta} (\vr,\vt) = \varrho(e_\delta(\varrho,\vt) - s_\delta (\varrho,\vt)) = H_1 (\vr,\vt) + \delta \vr (\vt - \log \vt)
%%$$
%\comment{\color{red} Here the choice of $\bar\vt$ in $H_{\bar\vt}$ simplifies our life now, but for example for the low Mach limit we probably can not expect $\bar\vt$ to be $1$, - but probably it doesn't matter, since we need to multyply entropy balance by $\bar\vt$ only.}

Since the vector field $\vc{V}$ is regular suitable manipulations with the H\"older, Young and Poincar|' e inequalities and thermodynamical hypothesis yield\footnote{Hereinafter, $c$ is a constant which is independent of solution. It depends on data, right hand side and it may vary from line to line. It may also depend on penalization parameters, however we emphasise this particular dependence as it play a role in further computations.}
\[
\int_B \vr \vu \cdot \vc{V} (\tau, \cdot)\, \dx  \leq c(\vc{V}) \int_B \sqrt\vr \sqrt \vr |\vu|(\tau, \cdot)\dx  \leq c + \frac 14 \int_B \vr |\vu|^2(\tau,\cdot)\dx,
\]
\begin{equation}\label{b_01}
\int_0^\tau \int_B  \mathbb S_\omega : \nabla_x \vc{V}\dxdt
 \leq c(\vc{V}, \lambda) + \int_0^\tau \int_B \frac{1}{6} \lambda \vt^5 \dxdt  + \frac 1 2 \int_0^\tau\int_B \frac{1}{\vt} \mathbb S_\omega:\Grad \vu \dxdt,
\end{equation}
%\comment{\color{red} here wee need to have $\int \int |\mathbb S|^2 \leq  \int \int \mathbb S : \nabla_x \vu$ !!! FOr the moment I don't see how to get it with $\eta \Div \vu$ term in $\mathbb S$ - Aneta;  in $C$ is hidden
%$\frac{1}{a_\eta^{1/3}}$ }
\[
\abs{\int_0^\tau \int_B \vr \vu \cdot \partial_t \vc{V}\dxdt}  \leq c \int_0^\tau \int_B \sqrt \vr \sqrt \vr |\vu|\, \dxdt
\leq c (\vc{V}, \vr_0)  + c\int_0^\tau \int_B \vr |\vu|^2 \, \dxdt ,
\]
\[
\abs{\int_0^\tau \int_B \vr (\vu\otimes \vu):\Grad \vc{V}\ {\rm d}x{\rm d}t}\leq  c\int_0^\tau \int_B \vr |\vu|^2 \, \dxdt .
\]
%\[
%\abs{\int_0^\tau\int_B \varrho \vc{f} \cdot (\vc{u} - \vc{V})\ {\rm d}x{\rm d}t}\leq c  + c \int_0^\tau \int_B \varrho|\vc{u}|^2\dxdt.
%\]
In order to deal with pressure term let us notice that
\begin{equation}\label{P>0}
P'(Z) > 0 \quad \mbox{ for all }Z\geq0.
\end{equation}
Indeed,  \eqref{pm_1} provides that $P'(Z) >0$ if $0<Z<\underline{Z}$ or $Z> \overline{Z}$. This together with \eqref{p2p} gives \eqref{P>0}. Next by \eqref{em1},  \eqref{pm2}, \eqref{pm3} one can infer that
	\begin{equation}\label{p_infty}
	\lim_{Z\to \infty} \frac{P(Z)}{Z^{\frac{5}{3}}} = p_\infty > 0.
	\end{equation}
Therefore by \eqref{p2p}, \eqref{pm2}, \eqref{pm3} we obtain the following bound on the molecular pressure
$p_M$
	$$\underline{c} \vr\vt \leq p_M \leq \overline{c} \vr\vt \quad \mbox{ if } \vr < \overline{Z} \vt^{\frac 3 2},$$
	\begin{equation}\label{pm_b}
	\underline{c} \vr^{\frac 5 3} \leq p_M \leq \overline{c}
\left\{
\begin{array}{ccc}
 \vt^{\frac 5 2} & \mbox{ if } \vr < \overline{Z} \vt^{\frac 3 2}  \\
\vr^{\frac 5 3}  &  \mbox{ if } \vr > \overline{Z} \vt^{\frac 3 2}  , \\
\end{array}
\right.
	\end{equation}
where we use also monotonicity of $p_M$ in $\vr$  to control it on the set $\underline{Z} \vt^{\frac 3 2} \leq \vr \leq \overline{Z} \vt^{\frac 3 2}$. With the above informations at hand we deduce
\begin{multline}\label{b_02}
\abs{\int_0^\tau \int_B p_{\eta,\delta}(\vr, \vt) \Div \vc{V} \dxdt}
 \leq
  c(\vc{V},\lambda) + c(\vc{V}) \int_0^\tau \int_B a_\eta \vt^4 \dxdt  + \int_0^\tau\int_B \frac{1}{6} \lambda \vt^5 \dxdt \\+ c(\vc{V})\int_0^\tau \int_B \vr^{\frac{5}{3}} \dxdt  + c(\vc{V}) \int_0^\tau \int_B \frac{\delta}{\beta-1} \varrho^\beta \dxdt.
\end{multline}
%\comment{In first $C$ is hiddent $\frac{1}{a_\eta^{5/3}}$}
%where $C(\vc{V})$ could also depend on $f$ and initial values of $\rho,\ \theta$ and $\vc{u}$.
%\comment{In order to use Gronwal inequality .... what left is to deal with  $\int_0^\tau \int_B p_{\delta}(\vr, \vt) \Div \vc{V} \dxdt $. }
Finally, as we may assume that $\lambda\leq 1$, we deduce by the Young inequality
$$
\int_{0}^\tau\int_B\lambda\vt^4 \dxdt\leq \lambda^{1/5} \int_0^\tau\int_B \lambda^{4/5}\vt^4 \dxdt \leq c + \frac16\int_0^\tau\int_B\lambda \vt^5 \dxdt .
$$

Moreover let us recall  that
 $$\vr e_\eta = \vr e_M + a_\eta \vt^4$$
and one can prove that
	\begin{equation}\label{e_est}
	\vr e_\eta \geq a_\eta \vt^4 + \frac{3 p_\infty}{2} \vr^{\frac 5 3} .
	\end{equation}
Indeed, by \eqref{pm2}, \eqref{pm3}, and \eqref{p_infty}	
	\begin{equation}\label{e_infty}
	\lim_{\vt \to 0^+} e_M(\vr,\vt) = \frac{3}{2} \vr^{\frac 2 3} p_\infty .
	\end{equation}
By \eqref{em1} $e_M$ is strictly increasing function of $\vt$ on $(0,\infty)$ for any fixed $\vr$. This together with \eqref{e1e}  and \eqref{e_infty} justify	that
	\begin{equation}\label{vr_e>}
	\vr e_\eta (\vr,\vt) \geq \frac{3 p_\infty}{2} \vr^{\frac 5 3} + a_\eta \vt^4
	\end{equation}
and consequently \eqref{e_est} holds. Furthermore, see \eqref{em2},
$e_M(\vr,\vt) = \underline{e}_M(\vr) + \int_0^\vt \frac{\partial{e_M}}{\partial \vt}(\vt,\tau) {\rm\,d}\tau$. This together with \eqref{em1} and  \eqref{e_infty} provides
	\bFormula{em_a}
	0\leq e_M(\vr,\vt) \leq c(\vr^{\frac{2}{3}} + \vt).
	\eF

%
%The term $\frac{\delta}{\vt^2}$ on the right hand side of inequality \eqref{ae1} can be absorbed by the term
%$\frac{\delta}{\vt^3}$ on the left hand side if $\vt < 1$, otherwise we may bound it by $\delta$ (or some constant).

Finally, summarising all above considerations by Gronwall inequality  we obtain the following
\begin{equation}\label{uni_est_1}
\begin{split}
& \int_{B}\left(\frac 12 \vr |\vu|^2 + H_ {1,\eta}(\vr,\vt) +   \frac{\delta}{\beta -1} \vr^\beta  \right)(\tau,\cdot)\ {\rm d}x + \frac 1\varepsilon \int_0^\tau \int_{\Gamma_t} |(\vu - \vc{V})\cdot \vc{n}|^2  \,{\rm d}S_x{\rm d}t
 \\ &
 + \int_{0}^\tau \int_B \frac 1\vt \left( \frac{1}{2}\mathbb S_\omega :\Grad \vu + \frac{\kappa_{\nu}(\vt)|\Grad \vt|^2}\vt \right) \dxdt + \int_{0}^\tau \int_B \frac{1}{2}\lambda \vt^5\dxdt
 \\ &
 \leq c(\vc{V},\lambda) \int_B \left(\frac 12 \frac{(\vr \vu)_{0,\delta}^2}{\vr_{0,\delta}}  + H_{1,\eta}(\vr_{0,\delta},\vt_{0,\delta}) +  \frac{\delta}{\beta -1} \vr_{0,\delta}^\beta  - (\vr \vu)_{0,\delta} \cdot \vc{V}(0,\cdot) + 1\right) \dx .
\end{split}
\end{equation}

Let us notice that the continuity equations provides that
\begin{equation}\label{aprior_1}
\int_{B} \vr(\tau)  = \int_{B} \vr(0) \mbox{ for all }t\in (0,T).
\end{equation}

Setting $\bar{\vr}$ constant such that  $\int_B (\vr - \bar\vr ) \dx = 0$ for a.a. $\tau \in [0,T)$ we may rewrite \eqref{uni_est_1}  as the following total dissipation inequality
\begin{multline}\label{uni_est_2}
 \int_{B}\left(\frac 12 \vr |\vu|^2 + H_{1,\eta}(\vr,\vt) - (\vr - \bar\vr) \frac{\partial H_{1,\eta}(\bar\vr,1)}{\partial \vr}
- H_{1,\eta}(\bar\vr,1)+   \frac{\delta}{\beta -1} \vr^\beta  \right)(\tau,\cdot)\ {\rm d}x + \frac 1\varepsilon \int_0^\tau \int_{\Gamma_t} |(\vu - \vc{V})\cdot \vc{n}|^2  \,{\rm d}S_x{\rm d}t
 \\
 + \int_{0}^\tau \int_B \frac 1\vt \left( \frac{1}{2}\mathbb S_\omega :\Grad \vu + \frac{\kappa_{\nu}(\vt)|\Grad \vt|^2}\vt  \right) \dxdt +  \int_0^\tau \int_B \frac{1}{2} \lambda \vt^5 \dxdt
 \\
 \leq c(\vc{V},\lambda) \int_B \left(\frac 12 \frac{(\vr \vu)_{0,\delta}^2}{\vr_{0,\delta}}  + H_{1,\eta}(\vr_{0,\delta},\vt_{0,\delta})
  +  \frac{\delta}{\beta -1} \vr_{0,\delta}^\beta  - (\vr \vu)_{0,\delta} \cdot \vc{V}(0,\cdot) + 1\right) \dx
  \\ - \int_B \left((\vr_{0,\delta} - \bar\vr) \frac{\partial H_{1,\eta}(\bar\vr,1)}{\partial \vr}
- H_{1,\eta}(\bar\vr,1)\right) \dx.
\end{multline}

In such form the left hand side of \eqref{uni_est_2} is nonnegative due to the hypothesis of thermodynamic stability \eqref{pm_1}, \eqref{em1}.
%{\color{red} we need to provide that $H_{1,0,\delta}$ is bounded}

%\comment{Below I believe that constants depends only on parameters expressed in parenthesis, please correct if I'm wrong ... - Aneta}
Directly from \eqref{uni_est_2} we obtain that

\begin{equation}\label{aprior_5}
\int_0^T\int_{\Gamma_t} |(\vu - \vc{V})\vc{n}|^2 {\rm \,dS}_{x} \dt \leq \varepsilon c(\lambda),
\end{equation}
\begin{equation}\label{aprior_7}
\mbox{ess}\sup_{\tau\in(0,T)} \|\delta \varrho^\beta(\tau,\cdot)\|_{L^{1}(B)}\leq c(\lambda),
\end{equation}
\begin{equation}\label{aprior_8}
\mbox{ess}\sup_{\tau\in(0,T)} \|\sqrt \vr \vu (\tau,\cdot) \|_{L^2(B)} \leq c(\lambda),
\end{equation}
\begin{equation}\label{aprior_11}
\left\|\lambda \vt^5 \right\|_{L^1((0,T)\times B)}  \leq c(\lambda).
\end{equation}
By \eqref{uni_est_2} we get also
\begin{equation*}\label{aprior_2_0}
\int_0^T\int_B \omega |\nabla_x \vu |^2 \dxdt \leq c(\lambda)
\end{equation*}
and then  \eqref{i4}, \eqref{mu_1}, the generalized Korn-Poincar\'e inequality (\cite[Proposition~2.1]{FeNo6}) and \eqref{V_0} give rise to
\begin{equation}\label{aprior_2}
\| \vu \|_{L^2((0,T)\times B)}  + \|\nabla_x \vu\|_{L^2((0,T)\times B)} \leq c(\lambda,\omega).
\end{equation}
Next by  \eqref{uni_est_2}, \eqref{kappa_2},  \eqref{kappa_3} we get
\begin{equation}\label{aprior_12}
\int_0^T \int_B \chi_\nu \left( | \nabla_x \log (\vt)|^2 + | \nabla_x \vt^{\frac 3 2}|^2 \right) \dxdt \leq c (\lambda) .
\end{equation}
Since $H_1$ is coercive and bounded from below by \cite[Proposition 3.2]{FeNo6} we get the following 	
%\comment{Maybe it is worth to recall here Proposition 3.2 ???}
\begin{equation}\label{aprior_9}
\mbox{ess}\sup_{\tau\in(0,T)} \|\vr e_\eta\|_{L^1(B)} \leq c(\lambda),
\end{equation}
and consequently
\begin{equation}\label{aprior_10}
\mbox{ess}\sup_{\tau\in(0,T)} \|a_\eta\vt^4(\tau,\cdot) \|_{L^1(B)}  \leq c(\lambda),
\end{equation}
\begin{equation}\label{aprior_6}
\mbox{ess}\sup_{\tau\in(0,T)} \|\varrho (\tau,\cdot) \|_{L^{\frac{5}{3}}(B)}\leq c(\lambda).
\end{equation}	
Then by \eqref{aprior_12}, \eqref{aprior_10}, and by Poincar\'e inequality (see \cite[Proposition~2.2]{FeNo6})
\begin{equation}\label{aprior_3}
\| \vt^\gamma \|_{L^2(0,T; W^{1,2} (B))} \leq c( \lambda, \nu) \quad \mbox{ for any } 1\leq \gamma \leq \frac{3}{2}.
\end{equation}
Moreover  by  \eqref{aprior_12}, \eqref{aprior_10}, \eqref{aprior_6}
\begin{equation}\label{aprior_4}
\| \log \vt \|_{L^2(0,T; W^{1,2} (B))}  \leq c( \lambda, \nu)
\end{equation}
(for more details see \cite[Chapter~2.2.4]{FeNo6}).
From \eqref{aprior_6} and \eqref{aprior_8} we deduce
\begin{equation}\label{p15_3}
\|\vr \vu \|_{L^\infty(0,T;L^{\frac 5 4 }(B))}\leq c(\lambda).
\end{equation}

Moreover, we use the technique based on the Bogovskii operator in order to derive the existence of $\pi>0$ fulfilling
\bFormula{p15}
\int\int_{K} \Big( p(\vr, \vt) \vr^\pi + \delta \vr^{\beta + \pi} \Big) \, \dxdt \leq c({K}),
\eF
for any compact $K \subset (0,T) \times {B}$ such that
\bFormula{p15_1}
{K} \cap \left( \cup_{ \tau \in [0,T] } \Big( \{ \tau \} \times \Gamma_\tau \Big) \right) = \emptyset.
\eF
It is worth pointing out that $\pi$ in \eqref{p15} can be chosen independently of $\varepsilon,\ \nu, \ \omega, \ \eta, \ \lambda$ and $\delta$.
For details we refer reader to \cite[Section 4.2]{EF71} or \cite{FP13}.

By hypothesis \eqref{pm_1}--\eqref{pm3} and Gibbs' relation one can deduce that
$$|  s_M(\vr,\vt) | \leq c (1+ | \log \vr| + | \log \vt |) \quad \mbox{ for all }\vr,\,\vt > 0 \mbox{ and some } c>0,$$
see  \cite[Section 3.2]{FeNo6} for details. Therefore, there exists $c>0$ such that
	\begin{equation}\label{p15_2}
	\vr s_\eta(\vr,\vt) \leq c (\vr + \vr|\log \vr| + \vr |\log \vt|+ a_\eta \vt^3) .
	\end{equation}
Relation \eqref{p15_2} together with \eqref{aprior_11},   \eqref{aprior_6},  \eqref{aprior_4} give rise to
	\begin{equation}\label{p15_4}
	\| \vr s_\eta (\vr,\vt) \|_{L^p((0,T)\times B)} \leq c( \lambda,\nu) \quad\mbox{ with some } p >1.
	\end{equation}
Moreover, the above combined with  \eqref{p15_3}, \eqref{aprior_2} and the Sobolev embedding theorem
yields
	\begin{equation}\label{p15_5}
	\|\vr s_\eta(\vr,\vt) \vu\|_{L^q((0,T)\times B)}
	\leq c( \lambda,\nu,\omega)
	\quad\mbox{ with some }q>1.
	\end{equation}
By \eqref{aprior_11}, \eqref{aprior_3} and by observation that (due to \eqref{kappa_2}, \eqref{kappa_3})
	$$\frac{\kappa_{\nu}(\vt)}{\vt} |\nabla_x \vt | \leq c \chi_\nu \left( | \nabla_x \log(\vt)| + \vt^{\frac{3}{2}} |\nabla_x \vt^{\frac{3}{2}} |\right),$$
we infer that
	\bFormula{p15_10}
	\left\| \frac{\kappa_\nu (\vt)}{\vt} \nabla_x \vt \right\|_{L^r((0,T) \times B)} \leq c( \lambda,\nu)
	\quad \mbox{ with certain }r>1.
	\eF
%Directly by  \eqref{uni_est_2}, \eqref{k_nu} we have
%	\bFormula{p15_11}
%	\delta \int_0^T\int_B \frac{\chi_\nu}{\vt^3} |\nabla_x \vt |^2 \dxdt \leq c(\eta) .
%	\eF
%Due to \eqref{aprior_10}, \eqref{p15_11},  \eqref{aprior_3} we infer that
%	\bFormula{p15_12}
%	\delta \int_0^T  \| \vt^{- \frac{1}{2}} (t, \cdot)\|^2_{W^{1,2}(B)} \dt \leq c (\eta).
%	\eF
Next, according to \eqref{em_a}, \eqref{aprior_6}, \eqref{aprior_11}
we have that
	\bFormula{p15_13_00}
	\| \vr e_\eta (\vr, \vt) \|_{L^1((0,T)\times B) } \leq c( \lambda)
	\eF
and due to \eqref{p15} higher integrability independent of $\delta$ can be obtained only away of $\Gamma_\tau$ interface, namely on sets as in \eqref{p15_1}.
%
%we obtain that
%	\bFormula{p15_13}
%	\| \vr e_\eta (\vr, \vt) \|_{L^p( K)} \leq c( \lambda) \quad \mbox{ with some }p>1 \mbox{ and }K \mbox{ as in \eqref{p15_1}} .
%	\eF

%
%\comment{The dependence of the constants on $\delta$ in the above estimates is due to dependence of initial data on this parameter in RHS of \eqref{uni_est_2} .... most of the cases. In cases when this dependence is crucial for $\delta \to 0$ it is somehow easy do deduce. Shall we skip this dependence in of $c$ on $\delta$?? }
%%%%%%%%%
%%%%%%%%%

\section{Singular limits}

\label{s}

In this section, we perform successively the singular limits  $\ep \to 0$, $\eta \to 0$, $\omega\to 0$, $\nu \to 0$, $\lambda \to 0$ and $\delta\to 0$.

\subsection{Penalization limit. Passing with $\ep \to 0$}\label{s:penlimit}
\subsubsection{Direct consequences of uniform bounds}\label{s:first_limits}

Firstly, we proceed with $\ep \to 0$ in \eqref{p3}, \eqref{p4}, \eqref{p5}, \eqref{p6} and \eqref{p7} while other parameters $\nu$, $\omega$, $\eta$, $\lambda$ and $\delta$ remain fixed.
Let $\{ \vre, \vue, \vte \}_{\ep > 0}$ be the corresponding sequence of weak solutions of the penalized problem given by Theorem \ref{Tm2}.

First of all, directly from \eqref{aprior_5}, we derive that
	\bFormula{s1}
	( \vu - \vc{V} ) \cdot \vc{n} (\tau, \cdot) |_{\Gamma_\tau} = 0 \ \mbox{for a.a.}\ \tau \in [0,T].
	\eF
in the limit as $\varepsilon \to 0$.

By \eqref{aprior_6} we have
	\begin{equation}\label{s1_01}
	\vr_\ep \to \vr \quad \mbox{ weakly-(*) in }L^\infty(0,T; L^{\frac{5}{3}}(B)) ,
	\end{equation}
 by \eqref{aprior_10}
	\begin{equation}\label{s1_02}
	\vt_\ep \to \vt \quad \mbox{ weakly-(*) in }L^\infty(0,T; L^{4}(B)),
	\end{equation}
and due to \eqref{aprior_3} we get
		\begin{equation}\label{s1_03}
	\vt_\ep \to \vt \quad \mbox{ weakly in }L^2(0,T; W^{1,2}(B)).
	\end{equation}
Due to \eqref{aprior_11}, \eqref{aprior_10}, and \eqref{aprior_3} we have also that
	\begin{equation}\label{s1_tem4}
	\vt_\ep^4 \to \overline{\vt^4} \quad
    \mbox{ weakly in } L^1((0,T)\times B),
	\end{equation}
    \begin{equation}\label{s1_tem5}
	\vt_\ep^5 \to \overline{\vt^5} \quad
    \mbox{ weakly in } L^1((0,T)\times B).
	\end{equation}
    Here and in the rest of the paper the bar denotes a weak limit of a composed or nonlinear function.

Then,  \eqref{aprior_6} together with the equation of continuity \eqref{p3}, imply that
	\begin{equation}\label{s1_0}
	\vre \to \vr \ \mbox{in}\ C_{\rm weak}([0,T] ; L^{\frac{5}{3}} (B)).
	\end{equation}
Next by \eqref{aprior_2}, up to a subsequence, we get
	\begin{equation}\label{s1_1}
	\vue \to \vu \ \mbox{weakly in}\ L^2(0,T; W^{1,2}_0 (B, \mathbb{R}^3)).
	\end{equation}
Then since $L^{\frac{5}{3}}(B) \hookrightarrow \hookrightarrow W^{-1,2}(B)$, by \eqref{s1_01}, \eqref{s1_1},
\eqref{s1_0}, \eqref{aprior_8}, we obtain
	\bFormula{s1a}
	\vre \vue  \to \vr \vu \ \mbox{weakly-(*) in}\ L^\infty(0,T; L^{\frac{5}{4}}(B;\mathbb{R}^3)).
	\eF
Due to the embedding $W^{1,2}_0 (B) \hookrightarrow L^6(B)$ we infer
	\bFormula{s1b}
	\vre \vue \otimes \vue \to \Ov{ \vr \vu \otimes \vu }
	\ \mbox{weakly in}\ L^2(0,T; L^{\frac{30}{29}} (B; 	\mathbb{R}^3)).
	\eF

From \eqref{p4}  we deduce that
	\bFormula{s1c}
	\vre \vue \to \vr \vu \ \mbox{in} \ C_{\rm weak}([T_1, T_2]; L^{\frac{5}{4}}(O;\mathbb{R}^3))
	\eF
for any space-time cylinder such that
	\bFormula{s1c2}
	(T_1, T_2) \times {O} \subset [0,T] \times B,\  [T_1, T_2] \times \Ov{O} \cap \cup_{\tau \in [0,T]}
	\left( \{ \tau \} \times \Gamma_\tau \right) =
	\emptyset.
	\eF
Since $L^{\frac{5}{4}}(B) \hookrightarrow\hookrightarrow W^{-1,2}(B)$, we conclude that
	\bFormula{s1d}
	\Ov{\vr \vu \otimes \vu} = \vr \vu \otimes \vu \ \mbox{a.a. in}\ (0,T) \times B.
	\eF
Next by \eqref{aprior_11}, \eqref{aprior_6}, \eqref{p15} and asymptotic behaviour of $p_M$ we obtian
	\begin{equation}\label{s1_e}
	p_{\eta,\delta} (\vr_\ep, \vt_\ep) = p_M(\vr_\ep, \vt_\ep) + \frac{a_\eta}{3} \vt^4_\ep + \delta \vr^{\beta}
	\to \overline{p_M(\vr, \vt)}  + \frac{a_\eta}{3} \overline{\vt^4}     + \delta \overline{\vr^\beta}
	\quad \mbox{ weakly  in } L^1(K) \mbox{ with }K \mbox{ as in \eqref{p15_1}}.
	\end{equation}
Due to \eqref{aprior_2}, \eqref{aprior_11}, \eqref{mu_1}
	\begin{equation}\label{s1_f}
	\tn{S}_\omega (\vt_\ep, \Grad \vu_\ep) \to \overline{\tn{S}_\omega (\vt, \Grad \vu)} \quad
	\mbox{ weakly in } L^{\frac{4}{3}} ((0,T) \times B).
	\end{equation}

According to \eqref{p15_4}
	\begin{equation}\label{s1_g}
	\vr_\varepsilon s_\eta(\vr_\varepsilon,\vt_\varepsilon)\to \overline{\vr s_\eta(\vr,\vt)}
	\quad\mbox{ weakly in }L^q((0,T) \times B) \mbox{ with some } q>1,
	\end{equation}
and by \eqref{p15_5}
	\begin{equation}\label{s1_h}
	\vr_\varepsilon s_\eta(\vr_\varepsilon,\vt_\varepsilon)\vu_\ep \to \overline{\vr s_\eta(\vr,\vt) \vu}
	\quad\mbox{ weakly in }L^p((0,T) \times B) \mbox{ with some } p>1.
	\end{equation}
%By \eqref{p15_13} we have
%		\begin{equation}\label{s1_i}
%		\vr_\ep e_\eta (\vr_\ep, \vt_\ep) \to \overline{\vr e_\eta (\vr, \vt)} \quad \mbox{ weakly in }
%		L^1((T_1,T_2)\times O)
%		\end{equation}
%for any space-time cylinder as in \eqref{s1c2}.

More details for the above considerations can be found  in \cite{EF70,FeNo6}.

\subsubsection{Pointwise convergence of the temperature and the density}\label{sec_ae_conv}

In order to show a.e. convergence of the temperature we follow \cite{FeNo6}. The proof is based on the Div-Curl Lemma (see \cite[Section 3.6.2]{FeNo6}) and Young measures methods. To this end we set
\begin{equation*}
\begin{split}
\vc{U}_\varepsilon& = \left[\vr_\varepsilon s_\eta(\vr_\varepsilon,\vt_\varepsilon), \vr_\varepsilon s_\eta(\vr_\varepsilon,\vt_\varepsilon)\vu_\varepsilon + \frac{\kappa_{\nu}(\vt_\varepsilon)\nabla\vt_\varepsilon}{\vt_\varepsilon}\right]\\
\vc{W}_\varepsilon& = [G(\vt_\varepsilon),0,0,0],
\end{split}
\end{equation*}
where $G$ is a bounded and Lipschitz function on $[0,\infty)$. Then due to estimates obtained in previous section ${\rm Div}_{t,x} \vc{U}_\ep$ is precompact in $W^{-1,s} ((0,T)\times B)$ and ${\rm Curl}_{t,x} \vc{W}_\ep$ is precompact in $W^{-1,s}((0,T)\times B)^{4\times 4}$ with certain $s>1$. Therefore we can deduce that
  assumptions  of Div-Curl Lemma for $\vc{U}_\varepsilon$ and $\vc{W}_\varepsilon$ are satisfied and  we may derive
	\bFormula{p15_6}
	\overline{\vr s_\eta (\vr,\vt) G(\vt)} = \overline{\vr s_\eta(\vr,\vt)}\ \overline{G(\vt)}.
	\eF
%Further,
%$$
%\vr s_\delta(\vr,\vt) = \vr s_M(\vr,\vt) + \frac 43 a \vt^3 + \delta\vr\log\vt.
%$$
Next step is to show that
	\begin{equation}\label{p15_7}
	\overline{\vr s_M(\vr,\vt) G(\vt)}\geq \overline{\vr s_M(\vr,\vt)}\ \overline{G(\vt)},
	\quad %%\overline{\vr \log(\vt) G(\vt)} \geq \overline{\vr \log(\vt)}  \overline{ G(\vt)} ,
	\quad \overline{\vt^3 G(\vt)}\geq \overline{\vt^3} \ \overline{G(\vt)}.
	\end{equation}
for any continuous and increasing function $G$.
It can be derived by application of the theory of parametrized Young measures. Details can be found in \cite[Section 3.6.2]{FeNo6}.
Combining \eqref{p15_6}, \eqref{p15_7} and taking  $G(\vt)=\vt$ we deduce
	$$
	\overline{\vt^4} = \overline{\vt^3}\ \vt
	$$
which yields
	\begin{equation}\label{temp.1}
	\vt_\varepsilon \rightarrow \vt \ \mbox{ a.a. in }(0,T)\times B.
	\end{equation}
Moreover due to \eqref{aprior_4} the limit temperature $\vt$ is positive a.e. on the set $(0,T)\times B$.
%\comment{\color{red} Here I will use a different argument - Aneta}
%{\color{red} Namely, by \eqref{temp.1}, \eqref{aprior_11} and weak lower semi-continuity of convex function
%	$$\vt^{-1} \in L^1((0,T)\times B).$$
%}
Similarly as in \cite[Section 4.1.2]{KrMaNeWr} we deduce
\bFormula{s8}
\vre \to \vr \ \mbox{a.e. in}\ (0,T) \times B.
\eF

%We are now ready to prove weak $L^1$ convergence of the pressure. Since pressure is bounded in $L^1$ we need to show that it is equi-integrable -- see \cite[Theorem 2.10]{EF70}. However, estimates \eqref{aprior_11} and \eqref{p15} yields equi-integrability locally in the interior of $Q_T$ and $Q_T^c$. In order to show the equi-integrability up to the boundary we adopt the method of \cite[Section 2.2.6]{EF71}. The presence of the moving boundary need some additional technical computation, nevertheless, it does not bring any new difficulties and thus we omit the proof. Thus, \eqref{temp.1} and \eqref{s8} yield
%\begin{equation}\label{press.weak}
%p_{\eta,\delta}(\vr_\varepsilon, \vt_\varepsilon)\rightarrow p_{\eta,\delta}(\vr,\vt),\ \mbox{weakly in}\ L^1((0,T)\times B).
%\end{equation}
%and similarly
%\begin{equation}\label{press.weak2}
%p_{\eta,\delta}(\vr_\varepsilon, \vt_\varepsilon)\rightarrow p_{\eta,\delta}(\vr,\vt),\ \mbox{weakly in}\ L^1((0,T)\times B\setminus\Omega_t)
%\end{equation}

%\comment{Isn't there any modification needed due to different form of a pressure function ? - Aneta}
\subsubsection{The limit system as  $\ep \to 0$}
Let us summarise our considerations from Section~\ref{s:first_limits},~\ref{sec_ae_conv}.
Passing to the limit in \eqref{p3} we obtain by \eqref{s1_0}, \eqref{s1a} that
\bFormula{s2}
\int_{B} \vr \varphi (\tau, \cdot) \, \dx - \int_{B} \vr_{0,\delta} \varphi (0, \cdot) \, \dx =
\int_0^\tau \int_{B} \left( \vr \partial_t \varphi + \vr \vu \cdot \Grad \varphi \right) \, \dxdt
\eF
for any $\tau \in [0,T]$ and any $\varphi \in \DC([0,T] \times \mathbb{R}^3)$. Moreover,  the limit solutions satisfies also the renormalized equation in the same for as \eqref{p3}.

Next we proceed to a limit in \eqref{p4}. Since we have at hand only the \emph{local estimates} on the pressure, see \eqref{p15}, \eqref{p15_1},
we have to restrict ourselves to the class of test functions
\bFormula{s4}
\vph \in C^1([0,T); W^{1, \infty}_0 (B; \mathbb{R}^3)),\ {\rm supp}[ \Div \vph (\tau, \cdot)] \cap \Gamma_\tau = \emptyset,\
\vph \cdot \vc{n}|_{\Gamma_\tau} = 0 \ \mbox{for all}\ \tau \in [0,T].
\eF

In accordance with  \eqref{s1_01}, \eqref{s1_tem4}, \eqref{s1_1}, \eqref{s1a}, \eqref{s1b}, \eqref{s1d}, \eqref{s1_e}, \eqref{s1_f}, \eqref{temp.1}, \eqref{s8}, assumptions on \eqref{mu_1}, \eqref{eta_1}, the momentum equation reads
\bFormula{s9}
\int_0^\tau \int_{B} \left( \vr \vu \cdot \partial_t \vph + \vr [\vu \otimes \vu] : \Grad \vph + {p_{\eta,\delta}(\vr,\vt)} \Div \vph
-  \mathbb S_\omega : \Grad \vph \right) \, \dxdt
\eF
\[
= %\int_0^T \int_B \vr {\bf f}\cdot \vph \, \dxdt
- \int_{B} (\vr \vu)_{0,\delta} \cdot \vph (0, \cdot) \, \dx
\]
for any test function $\vph$ as in (\ref{s4}).

Further, due to \eqref{aprior_3}, \eqref{p15_10}, and \eqref{temp.1} we get
$\frac{\kappa_\nu(\vt_\varepsilon)}{\vt_\varepsilon} |\nabla_x \vt_\ep| \to \frac{\kappa_\nu (\vt)}{\vt} | \nabla_x \vt|$ weakly in $L^1((0,T)\times B)$.
The terms $\frac{1}{\vt}\mathbb S_\omega(\vt_\ep,\nabla_x \vu_\ep):\Grad \vu_\ep$ and $\frac{\kappa_\nu(\vt_\ep) |\Grad \vt_\ep|^2}{\vt}$ for $\vt\geq 0$ are lower weakly semicontinuous.
Then this together with \eqref{s1_tem4}, \eqref{s1_tem5}, \eqref{s1_g}, \eqref{s1_h}, \eqref{temp.1},  \eqref{s8} allows to conclude that
	\begin{equation}\label{s10}
	\begin{split}
	\int_0^T\int_B  & \vr s_\eta(\vr,\vt) (\partial_t \varphi + \vu\cdot \Grad \varphi) \dxdt
	  - \int_0^T\int_B \frac{\kappa_\nu(\vt)\Grad\vt \cdot \Grad\varphi}{\vt}\dxdt
	 \\ &  + \int_0^T \int_B \frac{\varphi}{\vt}\left(\mathbb S_\omega : \Grad \vu
	+ \frac{\kappa_{\nu}(\vt)|\nabla_x \vt|^2}{\vt}\right) \dxdt    -  \int_0^T \int_B \lambda\vt^4 \dxdt
	\leq -\int_B (\vr s)_{0,\delta,\eta}\varphi(0)\dx
	\end{split}
	\end{equation}
as $\ep\to 0$.

Finally, we can proceed to a limit with $\ep \to 0$ in \eqref{p7}. Since the sequence $\{ \vr_\ep e_\eta(\vr_\ep, \vt_\ep) \}_\ep$ is nonnegative and \eqref{temp.1}, \eqref{s8}, \eqref{p15_13_00}, by the Fatou lemma we deduce
$$
\limsup\int_0^T\int_B \vr_\varepsilon e_\eta(\vr_\varepsilon,\vt_\varepsilon)\partial_t \psi \dxdt\leq \int_0^T\int_B \vr e_{\eta}(\vr,\vt)\partial_t \psi \dxdt
$$
as far as $\partial_t\psi\leq 0$. Using \eqref{Vselenoidal}, \eqref{s1_e} and arguments used also above we obtain
\begin{multline}\label{s11}
\int_0^T\int_B\left( \left(\frac 12 \vr |\vu|^2 + \vr e_\eta(\vr,\vt)  + \frac{\delta}{\beta-1} \vr^\beta \right) \partial_t \psi  -   \lambda \vt^5 \psi \right) \dxdt \\
\geq - \int_{B} \left(\frac 12 \frac{(\vr \vu)_0^2}{\vr_{0,\delta}} + \vr_{0,\delta} e_{0,\delta,\eta} + \frac{\delta}{\beta-1} \vr_{0,\delta}^\beta  - (\vr \vu)_0\cdot {\bf V}(0,\cdot)\right) \psi(0) \dx \\
-  \int_0^T\int_B\left(\mathbb S_{\omega} : \Grad \vc{V} \psi  - \vr \vu \cdot \partial_t( \vc{V} \psi) - \vr (\vu\otimes \vu):\Grad \vc{V}\psi  - p_{\eta,\delta}(\vr, \vt) \Div \vc{V}\psi  \right)\dxdt
\end{multline}
for all $\psi \in C^1_c([0,T))$, $\partial_t \psi\leq 0$.
%\comment{Not really sure how we proceed to a limit on the left hand side. Is the w.l.s.c. sufficient? -- Venca}

\subsection{Fundamental lemma and extending the class of test functions}
\label{a}
%\texttt{\color{red} (This fundamental lemma should came without changes: we have bounds on density in $L^2$ due to artificial pressure and it is based on the continuity equation only. Then starting with zero initial data for density outside of time-space ''cylinder'' we get density zero, and terms involving it vanishes on rigid part, in particular term involving pressure in momentum equation $p \Div \varphi$, because of the choice of constitutive relation on $p$.)\\}
In order to get rid of the density dependent terms supported by the "solid" part $\left((0,T)\times B\right)\setminus Q_T$ we use \cite[Lemma 4.1]{FKNNS} which reads as

\bLemma{a1}
Let $\vr \in L^\infty (0,T; L^2(B))$, $\vr \geq 0$,  $\vu \in L^2(0,T; W^{1,2}_0(B;\mathbb{R}^3))$ be a weak solution of the equation of continuity,
specifically,
\bFormula{a1}
\int_{B} \Big( \vr (\tau, \cdot) \varphi (\tau, \cdot) - \vr_0 \varphi(0, \cdot) \Big) \,\dx
= \int_0^\tau \int_{B} \Big( \vr \partial_t \varphi + \vr \vu \cdot \Grad \varphi \Big) \, \dx\dt
\eF
for any $\tau \in [0,T]$ and any test function $\varphi \in C^1_c ([0,T] \times \mathbb{R}^3)$.

In addition, assume that
\bFormula{a2}
( \vu - \vc{V} )(\tau, \cdot) \cdot \vc{n} |_{\Gamma_\tau}  = 0 \ \mbox{for a.a.}\ \tau \in (0,T),
\eF
and that
\[
\vr_0 \in L^2 (\mathbb{R}^3), \ \vr_0 \geq 0,  \ \vr_0 |_{B \setminus \Omega_0} = 0.
\]

Then
\[
\vr(\tau, \cdot) |_{B \setminus \Omega_\tau} = 0 \ \mbox{for any}\ \tau \in [0,T].
\]

\eL

\bigskip

Since we have set our initial density $\vr_{0,\delta}$  to be zero on $B \setminus \Omega_0$ (see \eqref{data1}), by virtue of Lemma \ref{La1},
the continuity equation \eqref{s2}
reads
\bFormula{s20_0}
\int_{\Omega_\tau} \vr \varphi (\tau, \cdot)  \dx - \int_{\Omega_0} \vr_{0,\delta} \varphi (0, \cdot)  \dx =
\int_0^\tau \int_{\Omega_t} \left( \vr \partial_t \varphi + \vr \vu \cdot \Grad \varphi \right) \dxdt
\eF
for any $\tau \in [0,T]$ and any $\varphi \in \DC([0,T] \times \mathbb{R}^3)$. Moreover the following renormalized formulation is satisfied
\bFormula{s20_00}
\int_0^T \int_{\Omega_t} \vr B(\vr) ( \partial_t \varphi + \vu \cdot \Grad \varphi )\, \dxdt =
\int_0^T \int_{\Omega_t} b(\vr)  \Div \vu \varphi\, \dxdt - \int_{\Omega_0}  \vr_{0,\delta} B(\vr_{0,\delta})  \varphi (0)\, \dx
\eF
 for any $\varphi \in C^1_c([0,T) \times \mathbb{R}^3)$, and any  $b \in L^\infty \cap C [0, \infty)$ such that { $b(0)=0$} and
 $B(\vr) = B(1) + \int_1^\vr \frac{b(z)}{z^2} {\rm d}z.$
 Next the momentum equation (\ref{s9}) reduces to
\bFormula{s20_1}
\int_0^T \int_{\Omega_t} \left( \vr \vu \cdot \partial_t \vph + \vr [\vu \otimes \vu] : \Grad \vph + {p_{\eta,\delta}(\vr,\vt)} \Div \vph
-  \mathbb S_\omega(\vt,\nabla_x \vu) : \Grad \vph \right)   \dxdt
\eF
\[
=  - \int_{\Omega_0} (\vr \vu)_{0,\delta} \cdot \vph (0, \cdot) \dx + \int_0^T \int_{B\setminus \Omega_t} \tn{S}_\omega (\vt,\nabla_x \vu):\nabla \vph \dxdt - \int_0^T \int_{B\setminus \Omega_t} \frac{a_\eta}{3} \vt^4\Div \vph \dxdt
\]
for any test function $\vph$ as in (\ref{s4}). We remark that in this step we crucially need the extra pressure term $\delta \vr^\beta$ ensuring the density $\vr$ to be square integrable (see \eqref{aprior_7}).

Next, we argue the same way as in \cite[Section 4.3.1]{FKNNS} to conclude, that the momentum equation \eqref{s20_1} holds in fact for any test function $\vph$ such that
\begin{equation}\label{s4a}
\vph \in C^\infty_c([0,T]\times B;\mathbb{R}^3), \qquad \vph(\tau, \cdot)\cdot \vc{n} |_{\Gamma_\tau} = 0 \quad \text{ for any } \tau \in [0,T].
\end{equation}
Moreover, by Lemma~\ref{La1} and the choice of initial data the balance of entropy \eqref{s10} takes the following form
	\begin{equation}\label{s20_2}
	\begin{split}
	\int_0^T\int_{\Omega_t}  & \vr s(\vr,\vt) (\partial_t \varphi + \vu\cdot \Grad \varphi) \dxdt
	+ \int_0^T\int_{B\setminus \Omega_t}
	\frac{4}{3} a_\eta \vt^3 (\partial_t \varphi + \vu\cdot \Grad \varphi) \dxdt
	\\&
	  - \int_0^T\int_{\Omega_t} \frac{\kappa_{\nu}(\vt)\Grad\vt \cdot \Grad\varphi}{\vt}\dxdt
	    - \int_0^T \int_{B\setminus \Omega_t}
	    \frac{\kappa_{\nu}(\vt)\Grad\vt \cdot \Grad\varphi}{\vt}\dxdt
	 \\ &  + \int_0^T \int_{\Omega_t} \frac{\varphi}{\vt}\left(\mathbb S_\omega : \Grad \vu
	+ \frac{\kappa_{\nu}(\vt)|\nabla_x \vt|^2}{\vt} \right) \dxdt
	 + \int_0^T \int_{B\setminus\Omega_t} \frac{\varphi}{\vt}\left(\mathbb S_\omega : \Grad \vu
	+ \frac{\kappa_{\nu}(\vt)|\nabla_x \vt|^2}{\vt} \right) \dxdt
	\\&
	-  \int_0^T \int_{B} \lambda\vt^4 \dxdt
	\leq -\int_{\Omega_0} (\vr s)_{0,\delta}\varphi(0)\dx - \int_{B \setminus\Omega_0} \frac{4}{3} a_\eta \vt_{0,\delta}^3\varphi(0)\dx
	\end{split}
	\end{equation}
for all $\varphi \in C^1_c ([0,T) \times \overline{B})$, $\varphi \geq 0$. Finally, total energy balance \eqref{s11}   reads
	\begin{equation}\label{s20_3}
	\begin{split}
	&\int_0^T \int_{\Omega_t} \left(\frac 12 \vr |\vu|^2 + \vr e(\vr,\vt)
	+ \frac{\delta}{\beta-1} \vr^\beta \right) \partial_t \psi \dxdt
	 + \int_0^T \int_{ B \setminus \Omega_t} a_\eta \vt^4 \partial_t \psi \dxdt
	 -  \int_0^T \int_{B} \lambda\vt^5 \psi \dxdt
 	\\ &\geq  - \int_{\Omega_0} \left(\frac 12 \frac{(\vr \vu)_{0,\delta}^2}{\vr_{0,\delta}}
	+ \vr_{0,\delta} e_{0,\delta} + \frac{\delta}{\beta-1} \vr_{0,\delta}^\beta
	- (\vr \vu)_{0,\delta}\cdot {\bf V}(0,\cdot)\right) \psi(0) \dx - \int_{B\setminus \Omega_0} a_\eta \vt_{0,\delta}^4\psi(0) \dx
 	\\ &
	-  \int_0^T\int_{\Omega_t} \left(
    \mathbb S_{\omega} : \Grad \vc{V} \psi
	- \vr \vu \cdot \partial_t (\vc{V} \psi ) - \vr (\vu\otimes \vu):\Grad \vc{V} \psi
	- p_{\delta}(\vr, \vt) \Div \vc{V} \psi
	\right)\dxdt
	\\ &
    -
 	\int_0^\tau\int_{B\setminus \Omega_t} \left(
 	\mathbb S_{\omega}: \Grad \vc{V}
 	- \frac{1}{3} a_\eta \vt^4 \Div \vc{V}
 	 \right) \psi \dxdt
	\end{split}
	\end{equation}
for all $\psi \in C^1_c([0,T))$, $\partial_t \psi \leq 0$.

\subsection{Limit in radiation $\eta \to 0$}\label{s:rad}

Let us denote by $\{ \vr_\eta, \vu_\eta,\vt_\eta\}_{\eta>0}$ solutions to the system \eqref{s20_0}, \eqref{s20_00},  \eqref{s20_1}, \eqref{s20_2}, \eqref{s20_3}.   In this section we pass to the limit with $\eta \to 0$ and get rid of radiative components of the pressure, internal entropy and internal energy functions outside of the fluid domain. Let us notice that estimates obtained in Section~\ref{s:bounds} are independent of parameter $\eta$ (if not emphasised). Therefore by \eqref{aprior_11} and as $a_\eta = \eta a$ on $B\setminus \Omega_\tau$ for $\tau \in [0,T]$
	\begin{equation}\label{lr_1}
	 \left| \int_0^T\int_{B\setminus \Omega_t} \frac{1}{3}a \eta\vt_\eta^4 \Div \vph \dxdt \right| \leq \eta c \| \vt_\eta^5 \|^{\frac{4}{5}}_{L^1((0,T)\times B)} \|\Div \vph  \|_{L^\infty((0,T)\times (B\setminus \Omega_t))} \to 0 \mbox{ as } \eta \to 0,
	\end{equation}
where $\vph$ is as in \eqref{s4a}.
In a similar way
	\begin{equation}\label{lr_2}
	\int_0^T\int_{B\setminus \Omega_t} \frac{1}{3} a_\eta\vt_\eta^4 \Div \vc{V} \dxdt \dxdt \to 0 \mbox{ as } \eta \to 0,
	\end{equation}
	\begin{equation}\label{lr_3}
	\int_0^T\int_{B\setminus \Omega_t} \frac{4}{3}a_\eta\vt_\eta^3 \partial_t \varphi \dxdt  \to 0 \mbox{ as } \eta \to 0 \quad \mbox{ for any }\varphi \in C^1_c ([0,T) \times \overline{B}).
	\end{equation}
Next by  \eqref{aprior_11},  \eqref{aprior_2}
	\begin{equation}\label{lr_4}
	\left| \int_0^T\int_{B\setminus \Omega_t} \frac{4}{3} a_\eta \vt_\eta^3 \vu_\eta\cdot\Grad \varphi \dxdt   \right|
	\leq \eta c(T,B) \| \vt^5_\eta \|^{\frac{3}{10}}_{{L^1((0,T)\times B)} }  \| \vu_\eta\|^{\frac{1}{2}}_{L^2(0,T; L^6(B))}\| \nabla_x \varphi \|_{L^\infty((0,T)\times (B\setminus \Omega_t))} \to 0 \mbox{ as } \eta \to 0
	\end{equation}
for any $\varphi \in C^1_c ([0,T) \times \overline{B})$.
%%(Now depending if we consider energy equality/inequality with or without test function $\psi$:
%%if with we can skip the term $$\int_{B\setminus \Omega_{t}} \eta a \vt^4 \dx (t)$$ but then we need to go for the inequality after tis passage, or we can pass here to zero if we keep test function using \eqref{aprior_10})
%%
Since $\vt_\eta \to \vt$ weakly in $L^1((0,T)\times B)$, due to \cite[Corollary~2.2]{EF70}, we obtain
	$$\int_0^T \int_{B} \lambda \vt^5 \dxdt
	\leq \lim\inf_{\eta\to 0} \int_0^T \int_{B} \lambda \vt^5_\eta \dxdt  .$$

The initial condition terms involving $a_\eta$ obviously converge to zero. To pass to the limit in remaining terms outside of the fluid part and in all terms in the fluid part we use the same arguments as for passing with  $\ep \to 0$. The continuity equation in the limit $\eta \to 0$ takes the same form as in \eqref{s20_0}, \eqref{s20_00}. The momentum equation \eqref{s20_1} as $\eta \to 0$ satisfies
\bFormula{s207_1}
\int_0^T \int_{\Omega_t} \left( \vr \vu \cdot \partial_t \vph + \vr [\vu \otimes \vu] : \Grad \vph + {p_\delta(\vr,\vt)} \Div \vph
-  \mathbb S_\omega(\vt,\nabla_x \vu) : \Grad \vph \right)   \dxdt
\eF
\[
=  - \int_{\Omega_0} (\vr \vu)_{0,\delta} \cdot \vph (0, \cdot) \, \dx - \int_0^T \int_{B\setminus \Omega_t} \tn{S}_\omega (\vt,\nabla_x \vu):\nabla \vph \dxdt
\]
for any test function $\vph$ as in \eqref{s4a}. The entropy inequality \eqref{s20_2} in the limit reads as
	\begin{equation}\label{s207_2}
	\begin{split}
	\int_0^T\int_{\Omega_t}  & \vr s(\vr,\vt) (\partial_t \varphi + \vu\cdot \Grad \varphi) \dxdt
	  - \int_0^T\int_{\Omega_t} \frac{\kappa_{\nu}(\vt)\Grad\vt \cdot \Grad\varphi}{\vt}\dxdt
	    - \int_0^T \int_{B\setminus \Omega_t}
	    \frac{\kappa_{\nu}(\vt)\Grad\vt \cdot \Grad\varphi}{\vt}\dxdt
	 \\ &  + \int_0^T \int_{\Omega_t} \frac{\varphi}{\vt}\left(\mathbb S_\omega : \Grad \vu
	+ \frac{\kappa_{\nu}(\vt)|\nabla_x \vt|^2}{\vt} \right) \dxdt
	 + \int_0^T \int_{B\setminus\Omega_t} \frac{\varphi}{\vt}\left(\mathbb S_\omega : \Grad \vu
	+ \frac{\kappa_{\nu}(\vt)|\nabla_x \vt|^2}{\vt} \right) \dxdt
	\\&
	-  \int_0^T \int_{B} \lambda\vt^4 \dxdt
	\leq -\int_{\Omega_0} (\vr s)_{0,\delta}\varphi(0)\dx
	\end{split}
	\end{equation}
for all $\varphi \in C^1_c ([0,T) \times \overline{B})$, $\varphi \geq 0$. Finally, total energy balance \eqref{s20_3} in the limit $\eta\to 0$ satisfies
	\begin{equation}\label{s207_3}
	\begin{split}
	&\int_0^T \int_{\Omega_t} \left(\frac 12 \vr |\vu|^2 + \vr e
	+ \frac{\delta}{\beta-1} \vr^\beta \right) \partial_t \psi \dxdt
	 -  \int_0^T \int_{B} \lambda\vt^5 \psi \dxdt
 	\\ & \geq  - \int_{\Omega_0} \left(\frac 12 \frac{(\vr \vu)_{0,\delta}^2}{\vr_{0,\delta}}
	+ \vr_{0,\delta} e_{0,\delta} + \frac{\delta}{\beta-1} \vr_{0,\delta}^\beta
	- (\vr \vu)_{0,\delta}\cdot {\bf V}(0,\cdot)\right) \psi(0) \dx
 	\\ &
	-  \int_0^T\int_{\Omega_t} \left(
    \mathbb S_{\omega} : \Grad \vc{V} \psi
	- \vr \vu \cdot \partial_t (\vc{V} \psi ) - \vr (\vu\otimes \vu):\Grad \vc{V} \psi
	- p_{\delta}(\vr, \vt) \Div \vc{V} \psi
	\right)\dxdt
	\\ &
    -
 	\int_0^\tau\int_{B\setminus \Omega_t}
 	\mathbb S_{\omega} (\vt,\nabla_x \vu): \Grad \vc{V}
 	 \psi \dxdt
	\end{split}
	\end{equation}
for all $\psi \in C^1_c([0,T))$, $\psi \geq 0$, $\partial_t \psi \leq 0$.

\subsection{Vanishing viscosity $\omega \to 0$}\label{s:viscosity}
Let us denote by $\{ \vr_\omega, \vu_\omega, \vt_\omega \}_{\omega>0}$ solutions to
\eqref{s20_0}, \eqref{s20_00},
\eqref{s207_1}, \eqref{s207_2}, and \eqref{s207_3}. Now our aim is to pass with $\omega \to 0$ in order to get rid of terms related to the viscous stress tensor outside of the fluid domain.
Let us notice that again we may obtain analogous estimates as in Section~\ref{s:bounds} which are independent of $\omega$  if not emphasised.
For the viscous term in the momentum equation \eqref{s207_1} we observe that by \eqref{uni_est_2}, \eqref{aprior_11}
for any $\vph$ as in \eqref{s4a} we have
\begin{equation}\label{svv_1}
\begin{split}
\int_0^T & \int_{B\setminus \Omega_t}  \mathbb S_\omega(\vt_\omega,\nabla_x \vu_\omega) :\Grad \vph \dxdt \leq \sqrt{\omega}\int_0^T\int_{B\setminus \Omega_t} \frac 1{\sqrt{\vt_\omega} }
\sqrt{\omega} {\mathbb S}(\vt_\omega,\nabla_x \vu_\omega) \sqrt{\vt_\omega} :\Grad \vph \dxdt
\\ &
\leq \sqrt{\omega}\frac{1}{2} \int_0^T \int_{B\setminus \Omega_t} \frac{1}{\vt_\omega} \omega | \mathbb S(\vt_\omega, \nabla_x \vu_\omega)|^2 \dxdt +  \sqrt{\omega} \frac{1}{2} \int_0^T \int_{B\setminus \Omega_t} \vt_\omega |\Grad \vph|^2 \dxdt
\\ &
\leq \sqrt{\omega} c \int_0^T \int_{B\setminus \Omega_t} \frac{1}{\vt_\omega}  \mathbb S_\omega(\vt_\omega,\nabla_x \vu_\omega) :\nabla_x \vu_\omega  \dxdt  + \sqrt{\omega} c(\vph,\lambda) \int_0^T \int_{B\setminus \Omega_t} \lambda \vt_\omega^5 \dxdt  \leq \sqrt{\omega} c \to 0 \mbox{ as } \omega\to 0.
\end{split}
\end{equation}
In a similar way we show that in the total energy inequality \eqref{s207_3}
\begin{equation}\label{svv_10}
\int_0^T \int_{B\setminus \Omega_t} \mathbb S_{\omega} (\vt_\omega,\nabla_x \vu_\omega): \Grad \vc{V}   \dxdt
\to 0 \quad \mbox{ as }\omega\to 0.
\end{equation}
When passing with $\omega\to 0$ in entropy inequality \eqref{s207_2} we skip the term
$\int_0^T \int_{B\setminus\Omega_t} \frac{\varphi}{\vt} \mathbb S_\omega : \Grad \vu \dxdt$ since it is positive for all $\varphi \in C^1_c ([0,T) \times \overline{B})$, $\varphi \geq 0$.
Since on the fluid domain all estimates obtained in Section~\ref{s:bounds} holds true, we can pass tho the limit with $\omega \to 0$ in all terms on the fluid domain and on the remaining one outside of the fluid domain in the same way as in previous Section~\ref{s:rad}.

All the above arguments allow as to pass with $\omega \to 0$. Then the continuity equation takes the same form as in \eqref{s20_0}, \eqref{s20_00}.
The momentum equation   \eqref{s207_1} in the limit $\omega \to 0$ reads
\bFormula{s208_1}
\int_0^T \int_{\Omega_t} \left( \vr \vu \cdot \partial_t \vph + \vr [\vu \otimes \vu] : \Grad \vph + {p_\delta(\vr,\vt)} \Div \vph
-  \mathbb S(\vt,\nabla_x \vu) : \Grad \vph \right)   \dxdt
=  - \int_{\Omega_0} (\vr \vu)_{0,\delta} \cdot \vph (0, \cdot) \dx
\eF
for any test function $\vph$ as in \eqref{s4a}.
The entropy inequality \eqref{s207_2} as $\omega \to 0$ satisfies
	\begin{equation}\label{s208_2}
	\begin{split}
	\int_0^T\int_{\Omega_t}  & \vr s(\vr,\vt) (\partial_t \varphi + \vu\cdot \Grad \varphi) \dxdt
	  - \int_0^T\int_{\Omega_t} \frac{\kappa_{\nu}(\vt)\Grad\vt \cdot \Grad\varphi}{\vt}\dxdt
	    - \int_0^T \int_{B\setminus \Omega_t}
	    \frac{\kappa_{\nu}(\vt)\Grad\vt \cdot \Grad\varphi}{\vt}\dxdt
	 \\ &  + \int_0^T \int_{\Omega_t} \frac{\varphi}{\vt}\left(\mathbb S : \Grad \vu
	+ \frac{\kappa_{\nu}(\vt)|\nabla_x \vt|^2}{\vt} \right) \dxdt
	 + \int_0^T \int_{B\setminus\Omega_t} \frac{\varphi}{\vt}\left(\frac{\kappa_{\nu}(\vt)|\nabla_x \vt|^2}{\vt} \right) \dxdt
	\\&
	-  \int_0^T \int_{B} \lambda\vt^4 \dxdt
	\leq -\int_{\Omega_0} (\vr s)_{0,\delta}\varphi(0)\dx
	\end{split}
	\end{equation}
for all $\varphi \in C^1_c ([0,T) \times \overline{B})$, $\varphi \geq 0$. Finally, for the total energy balance \eqref{s207_3} we get for $\omega \to 0$ that
	\begin{equation}\label{s208_3}
	\begin{split}
	&\int_0^T \int_{\Omega_\tau} \left(\frac 12 \vr |\vu|^2 + \vr e
	+ \frac{\delta}{\beta-1} \vr^\beta \right) \partial_t \psi \dxdt
	 -  \int_0^T \int_{B} \lambda\vt^5 \psi \dxdt
 	\\ & \geq  - \int_{\Omega_0} \left(\frac 12 \frac{(\vr \vu)_{0,\delta}^2}{\vr_{0,\delta}}
	+ \vr_{0,\delta} e_{0,\delta} + \frac{\delta}{\beta-1} \vr_{0,\delta}^\beta
	- (\vr \vu)_{0,\delta}\cdot {\bf V}(0,\cdot)\right) \psi(0) \dx
 	\\ &
	-  \int_0^T\int_{\Omega_t} \left(
    \mathbb S : \Grad \vc{V} \psi
	- \vr \vu \cdot \partial_t (\vc{V} \psi ) - \vr (\vu\otimes \vu):\Grad \vc{V} \psi
	- p_{\delta}(\vr, \vt) \Div \vc{V} \psi
	\right)\dxdt
	\end{split}
	\end{equation}
for all $\psi \in C^1_c([0,T))$, $\psi \geq 0$, $\partial_t \psi \leq 0$.

\subsection{Vanishing conductivity $\nu \to 0$}\label{s:cond}
Let $\{ \vr_\nu, \vu_\nu, \vt_\nu \}_{\nu>0}$ denote solutions to the system obtained in the previous Section~\ref{s:viscosity},  namely satisfying \eqref{s20_0}, \eqref{s20_00}, \eqref{s208_1}, \eqref{s208_2}, and \eqref{s208_3} for each fixed $\nu >0$.
In this section we pass with $\nu \to 0$
what allow us to show that  the term
$\int_0^T\int_{B \setminus \Omega_t} \frac{\kappa_\nu(\vt_\nu)\Grad \vt}{\vt} \cdot \Grad \varphi \dxdt$ vanishes in the limit for any $\varphi \in C^1_c ([0,T) \times \overline{B})$.

Let us notice that for each fixed $\nu >0$ we still keep that
	$$
	\left\| \frac{\chi_\nu \kappa(\vt_\nu)|\Grad \vt_\nu|^2}{\vt_\nu^2}\right\|_{L^1((0,T) \times B)}
	\leq 	c(\lambda) \quad \mbox{ and } \quad
	\| \lambda \vt_\nu^5 \|_{L^1((0,T) \times B)} \leq c(\lambda)
	$$
independent w.r.t. $\nu$. Therefore due to \eqref{kappa_1}--\eqref{kappa_3} we deduce that
	\[
	\begin{split}
	\int_0^T\int_{B\setminus \Omega_t} \frac{\kappa_\nu(\vt_\nu)\Grad \vt_\nu}{\vt_\nu} \cdot\Grad \varphi
	\dxdt
	& \leq \sqrt{\nu}\left(\int_0^T\int_{B\setminus\Omega_t}
	\frac{\nu \kappa(\vt_\nu)|\Grad \vt_\nu|^2}{\vt_\nu^2} \dxdt \right)^{\frac{1}{2}}
	\left(\int_0^T\int_{B\setminus \Omega_t} |\Grad\varphi|^2 \kappa(\vt_\nu) \dxdt \right)^\frac 12
	\\ & \to 0 \quad \mbox{ as } \nu \to 0 .
	\end{split}
	\]
Since $\int_0^T \int_{B\setminus\Omega_t} \frac{\varphi}{\vt_\nu}\left(\frac{\kappa_{\nu}(\vt_\nu)|\nabla_x \vt_\nu|^2}{\vt_\nu} \right) \dxdt $ is positive for any $\varphi \in C^1_c ([0,T) \times \overline{B})$, $\varphi \geq 0$  and any $\nu>0$, we can skip this term in the vanishing conductivity limit of the internal entropy inequality.
Next proceeding in the analogous way as in Section~\ref{s:rad} we may pass with $\nu \to 0$ in the  remaining terms of \eqref{s20_0}, \eqref{s20_00}, \eqref{s208_1}, \eqref{s208_2}, \eqref{s208_3}. Therefore  we obtain that the continuity equation in the limit satisfies again  \eqref{s20_0}, \eqref{s20_00}.
The momentum equation takes the same for as in  \eqref{s208_1} as $\nu \to 0$.
The entropy inequality \eqref{s208_2} in the limit $\nu \to 0$ reads
	\begin{equation}\label{s209_2}
	\begin{split}
	\int_0^T\int_{\Omega_t}  & \vr s(\vr,\vt) (\partial_t \varphi + \vu\cdot \Grad \varphi) \dxdt
	  - \int_0^T\int_{\Omega_t} \frac{\kappa(\vt)\Grad\vt \cdot \Grad\varphi}{\vt}\dxdt
	 \\ &  + \int_0^T \int_{\Omega_t} \frac{\varphi}{\vt}\left(\mathbb S : \Grad \vu
	+ \frac{\kappa(\vt)|\nabla_x \vt|^2}{\vt} \right) \dxdt
	-  \int_0^T \int_{B} \lambda\vt^4 \dxdt
	\leq -\int_{\Omega_0} (\vr s)_{0,\delta}\varphi(0)\dx
	\end{split}
	\end{equation}
for all $\varphi \in C^1_c ([0,T) \times \overline{B})$, $\varphi \geq 0$.
The energy inequality in the limit $\nu \to 0$ have the same form as in  \eqref{s208_3}.

\subsection{Vanishing additional temperature term, $\lambda \to 0$}\label{s:lambda}

In this section we get rid the term related to coefficient $\lambda$ - the only terms which are left also outside of a fluid domain. Let $\{ \vr_\lambda, \vu_\lambda , \vt_\lambda\}_{\lambda \in(0,1)}$ be solution to the limit system obtained in previous Section~\ref{s:cond}. In order to pass with $\lambda \to 0$ in \eqref{s20_0}, \eqref{s20_00}, \eqref{s208_1}, \eqref{s209_2}, and \eqref{s208_3} we need to provide uniform estimates analogous to these obtained in Section~\ref{s:bounds}, but independent of $\lambda$. To this end we proceed in a similar way as therein, we only need to modify \eqref{b_01} and \eqref{b_02} as follows
(for convenience we skip the subscript $\lambda$ in below notation)
\[
\int_0^\tau \int_{\Omega_t}  \mathbb S(\vt, \nabla_x \vu) : \nabla_x \vc{V}\dxdt
 \leq c(\vc{V}) + \int_0^\tau \int_{\Omega_t} a\vt^4 \dxdt  + \frac 1 2 \int_0^\tau\int_{\Omega_t} \frac{1}{\vt} \mathbb S:\Grad \vu \dxdt,
\]
\begin{multline*}
\abs{\int_0^\tau \int_{\Omega_t} p_{\delta}(\vr, \vt) \Div \vc{V} \dxdt}
 \leq
  c(\vc{V}) + c(\vc{V}) \int_0^\tau \int_{\Omega_t} a \vt^4 \dxdt + c(\vc{V})\int_0^\tau \int_{\Omega_t} \vr^{\frac{5}{3}} \dxdt  + c(\vc{V}) \int_0^\tau \int_{\Omega_t} \frac{\delta}{\beta-1} \varrho^\beta \dxdt.
\end{multline*}
Consequently we may easy obtain that
 \begin{multline}\label{uni_est_lambda_2}
 \int_{\Omega_\tau}\left(\frac 12 \vr |\vu|^2 + H_1(\vr,\vt) - (\vr - \bar\vr) \frac{\partial H_1(\bar\vr,1)}{\partial \vr}
- H_1(\bar\vr,1)+   \frac{\delta}{\beta -1} \vr^\beta  \right)(\tau,\cdot)\ {\rm d}x
 \\
 + \int_{0}^\tau \int_{\Omega_t} \frac 1\vt \left( \frac{1}{2}\mathbb S :\Grad \vu + \frac{\kappa(\vt)|\Grad \vt|^2}\vt  \right) \dxdt +  \int_0^\tau \int_{B} \frac{1}{2} \lambda \vt^5 \dxdt
 \\
 \leq c(\vc{V}) \int_{\Omega_0} \left(\frac 12 \frac{(\vr \vu)_{0,\delta}^2}{\vr_{0,\delta}}  + H_1(\vr_{0,\delta},\vt_{0,\delta})
  +  \frac{\delta}{\beta -1} \vr_{0,\delta}^\beta  - (\vr \vu)_{0,\delta} \cdot \vc{V}(0,\cdot) + 1\right) \dx \\
  - \int_{\Omega_0} \left((\vr_{0,\delta} - \bar\vr) \frac{\partial H_1(\bar\vr,1)}{\partial \vr}
- H_1(\bar\vr,1)\right) \dx.
\end{multline}
Hence the following estimates hold
\begin{equation}\label{aprior_7_bis}
\mbox{ess}\sup_{\tau\in(0,T)} \|\delta \varrho^\beta (\tau,\cdot) \|_{L^{1}(\Omega_\tau)}
+
\mbox{ess}\sup_{\tau\in(0,T)} \|\sqrt \vr \vu(\tau,\cdot)\|_{L^2(\Omega_\tau)} \leq c,
\end{equation}
\begin{equation}\label{aprior_11_bis}
\left\| \lambda \vt^5 \right\|_{L^1((0,T)\times B)}  \leq c,
\end{equation}
\begin{equation}\label{aprior_2_bis}
\| \vu \|_{L^2(Q_T)}  + \|\nabla_x \vu\|_{L^2(Q_T)} \leq c,
\end{equation}
\begin{equation}\label{aprior_10_bis}
\mbox{ess}\sup_{\tau\in(0,T)} \|a\vt^4  (\tau,\cdot) \|_{L^1(\Omega_\tau)}
+
\mbox{ess}\sup_{\tau\in(0,T)} \|\varrho  (\tau,\cdot) \|_{L^{\frac{5}{3}}(\Omega_\tau)}\leq c.
\end{equation}	
\begin{equation}\label{aprior_12_bis}
\int_0^T \int_{\Omega_t}  \left( | \nabla_x \log (\vt)|^2 + | \nabla_x \vt^{\frac 3 2}|^2 \right) \dxdt \leq c,
\end{equation}
\bFormula{p15_bis}
\int\int_{K} \Big( p(\vr, \vt) \vr^\pi + \delta \vr^{\beta + \pi} \Big) \, \dxdt \leq c({K})
 \mbox{ for certain }\pi>0
\eF
\begin{equation}\label{p15_2_bis}
\mbox{and for any compact }K \subset Q_T \mbox{ such that }
{K} \cap \left( \cup_{ \tau \in [0,T] } \Big( \{ \tau \} \times \Gamma_\tau \Big) \right) = \emptyset,
\end{equation}
and
	\begin{equation}\label{p15_4_bis}
	\| \vr s (\vr,\vt) \|_{L^q(Q_T)} +
	\|\vr s(\vr,\vt) \vu\|_{L^q(Q_T)}
	+
	\left\| \frac{\kappa_\nu (\vt)}{\vt} \nabla_x \vt \right\|_{L^q(Q_T)} \leq c
	\quad \mbox{ with certain }q>1,
	\eF
	\bFormula{p15_13_bis}
	\| \vr e (\vr, \vt) \|_{L^1( Q_T)} \leq c \quad \mbox{ and } \quad \| \vr e (\vr, \vt) \|_{L^p( K)} \leq c \quad \mbox{ with some }p>1 \mbox{ and }K \mbox{ as in \eqref{p15_2_bis}}.
	\eF
Then due to \eqref{aprior_11_bis}
	$$\int_0^T \int_B \lambda \vt^4_\lambda \dxdt \to 0 \quad \mbox{ as } \lambda\to 0.$$
Next notice that  the term $\int_0^\tau \int_B \lambda \vt^5_\lambda\psi  \dxdt$ is non-negative for all $\lambda>0$ and all $\psi \in C^1_c([0,T))$, $\psi \geq 0$
and therefore can be skipped in the total energy inequality \eqref{s208_3}.
In all other terms we pass in the same way as in Section ~\ref{s:penlimit}. Not that \eqref{aprior_10_bis} provides enough informations in steps where \eqref{aprior_11} is used therein. Consequently in the limit
$\lambda \to 0$ the continuity equation and the momentum equation take the same fore as in
 \eqref{s20_0}, \eqref{s20_00}, \eqref{s208_1}, respectevly. The entropy inequality satisfies
	\begin{equation}\label{s2010_2}
	\begin{split}
	\int_0^T\int_{\Omega_t}  & \vr s(\vr,\vt) (\partial_t \varphi + \vu\cdot \Grad \varphi) \dxdt
	  - \int_0^T\int_{\Omega_t} \frac{\kappa(\vt)\Grad\vt \cdot \Grad\varphi}{\vt}\dxdt
	   + \int_0^T \int_{\Omega_t} \frac{\varphi}{\vt}\left(\mathbb S : \Grad \vu
	+ \frac{\kappa(\vt)|\nabla_x \vt|^2}{\vt} \right) \dxdt
	\\ & \leq -\int_{\Omega_0} (\vr s)_{0,\delta}\varphi(0)\dx
	\end{split}
	\end{equation}
for all $\varphi \in C^1_c ([0,T) \times \overline{B})$, $\varphi \geq 0$ as $\lambda \to 0$.
And for the total energy inequality we get for $\lambda \to 0$ that
	\begin{equation}\label{s2010_3}
	\begin{split}
	\int_0^T \int_{\Omega_\tau} \left(\frac 12 \vr |\vu|^2 + \vr e
	+ \frac{\delta}{\beta-1} \vr^\beta \right) &  \partial_t \psi \dxdt
 	 \geq  - \int_{\Omega_0} \left(\frac 12 \frac{(\vr \vu)_{0,\delta}^2}{\vr_{0,\delta}}
	+ \vr_{0,\delta} e_{0,\delta} + \frac{\delta}{\beta-1} \vr_{0,\delta}^\beta
	- (\vr \vu)_{0,\delta}\cdot {\bf V}(0,\cdot)\right) \psi(0) \dx
 	\\ &
	-  \int_0^T\int_{\Omega_t} \left(
    \mathbb S : \Grad \vc{V} \psi
	- \vr \vu \cdot \partial_t (\vc{V} \psi ) - \vr (\vu\otimes \vu):\Grad \vc{V} \psi
	- p_{\delta}(\vr, \vt) \Div \vc{V} \psi
	\right)\dxdt
	\end{split}
	\end{equation}
for all $\psi \in C^1_c([0,T))$, $\psi \geq 0$, $\partial_t \psi \leq 0$.

\subsection{Conclusion of the proof -- artificial pressure and temperature term}

We proceed with  $\delta$ to $0$ similarly as in \cite{FeNo6}. Finally in the energy inequality we choose a test function in a form of \eqref{psi_zeta} and pass with $\xi \to 0$. Consequently  obtain that our system satisfies the form
required in the Theorem~\ref{Tm1}.

\section{Discussion}
\label{d}
\begin{itemize}
\item Let is point out  that we considered the full slip boundary condition for the velocity field. We can consider the Dirichlet condition for the velocity field. In that case it will be examined by means of Birkman's penalized method. For more details see \cite{{FeNeSt}}.
\item We can also consider the Navier type of boundary conditions.
\item The regularity of domain $\Omega$ follows from applications of level set method in  the Fundamental Lemma \ref{La1}.
\item The essential in our paper was considering energy inequality instead of energy equation together with
 introducing the term $\lambda \vt ^5$ in the energy balance and the term $\lambda \vt ^4$ into the entropy balance.
\end{itemize}

\begin{Remark}\label{divboundry}
The condition \eqref{Vselenoidal} is not restrictive. Indeed, for a general $\vc{V}\in C^1([0,T]; C_c^3(\mathbb R^3,\mathbb R^3))$ we can find $\vc{w}\in W^{1,\infty}(Q_T)$ such that $(\vc{V}-\vc{w})|_{\Gamma_\tau} = 0$ for all $\tau \in [0,T]$ such that $\Div \vc{w} = 0$ on some neighborhood of $\Gamma_\tau$ -- see \cite[Section 4.3.1]{FKNNS}. This function can be used in place of $\vc{V}$ in the definition of weak solutions and later also in the approximate problem. Due to the fact that $\vph = \vc{V} - \vc{w}$ is a suitable test function in \eqref{m3} and \eqref{p4}, the appropriate energy balances remain valid.
\end{Remark}

\def\ocirc#1{\ifmmode\setbox0=\hbox{$#1$}\dimen0=\ht0 \advance\dimen0
  by1pt\rlap{\hbox to\wd0{\hss\raise\dimen0
  \hbox{\hskip.2em$\scriptscriptstyle\citrc$}\hss}}#1\else {\accent"17 #1}\fi}

%\bibliography{citace}

\begin{thebibliography}{10}

\bibitem{AW}
 S.~ S. Antman, J.~P.~ Wilber.
 \newblock{The asymptotic problem for the springlike motion of a heavy piston in a viscous gas.}
 {\em Quart. Appl. Math.},{\bf 65}, 3: 471–-498, 2007.

\bibitem{BrKrMa}
J.~B\v rezina, O.~Kreml, and V.~M\'acha
\newblock Dimension reduction for the full Navier-Stokes-Fourier system
\newblock {\em J. Math. Fluid Mech.}, {\bf 19}:659--683, 2017.

\bibitem{BMR1}
M.~Bul{\' \i}{\v c}ek, J.~M{\' a}lek, and K.R. Rajagopal.
\newblock Navier's slip and evolutionary {N}avier-{S}tokes-like systems with
  pressure and shear-rate dependent viscosity.
\newblock {\em Indiana Univ. Math. J.}, {\bf 56}:51--86, 2007.

%\bibitem{Cor}
%F.~Coron.
%\newblock Derivation of slip boundary conditions for the {N}avier-{S}tokes
%  system from the {B}oltzmann equation.
%\newblock {\em J. Statistical Phys.}, {\bf 54}:829--857, 1989.

\bibitem{DL}
R.J. DiPerna and P.-L. Lions.
\newblock Ordinary differential equations, transport theory and {S}obolev
  spaces.
\newblock {\em Invent. Math.}, {\bf 98}:511--547, 1989.

%\bibitem{EF61}
%E.~Feireisl.
%\newblock Compressible {N}avier-{S}tokes equations with a non-monotone
%  pressure law.
%\newblock {\em J. Differential Equations}, {\bf 184}:97--108, 2002.

\bibitem{EF70}
E.~Feireisl.
\newblock {\em Dynamics of viscous compressible fluids}.
\newblock Oxford University Press, Oxford, 2004.

\bibitem{EF71}
E.~Feireisl.
\newblock On the motion of a viscous, compressible, and heat conducting fluid.
\newblock {\em Indiana Univ. Math. J.}, {\bf 53}:1707--1740, 2004.

\bibitem{BumJa}
E. Feireisl, B. J. Jin, and A. Novotn\'y.
\newblock Relative entropies, suitable weak solutions, and weak-strong uniqueness for the compressible Navier-Stokes system.
\newblock {\em J. Math. Fluid Mech. 14}, no. 4, 717--730, 2012.

\bibitem{FeNeSt}
E.~Feireisl, J.~Neustupa, and J.~Stebel.
\newblock Convergence of a {B}rinkman-type penalization for compressible fluid
  flows.
\newblock {\em J. Differential Equations}, {\bf 250}(1):596--606, 2011.

\bibitem{FKNNS} E. Feireisl, O. Kreml, \v{S}. Ne\v casov\'a, J. Neustupa, J. Stebel. {\em Weak solutions to the barotropic Navier-Stokes system with slip boundary conditions in time dependent domains.} J. Differential Equations \textbf{254} no. 1, 125--140, 2013.

\bibitem{FeNo6}
E.~Feireisl and A.~Novotn{\' y}.
\newblock {\em Singular limits in thermodynamics of viscous fluids}.
\newblock Birkh{\" a}user-Verlag, Basel, 2009.

\bibitem{FeNows}
E.~Feireisl and A.~Novotn\'y
\newblock Weak-strong uniqueness property for the full Navier-Stokes-Fourier System
\newblock {\em Arch. Rat. Mech. Anal.}, {\bf 204}:683--706, 2012.

\bibitem{FNP}
E.~Feireisl, A.~Novotn{\' y}, and H.~Petzeltov{\' a}.
\newblock On the existence of globally defined weak solutions to the
  {N}avier-{S}tokes equations of compressible isentropic fluids.
\newblock {\em J. Math. Fluid Mech.}, {\bf 3}:358--392, 2001.

\bibitem{FP13}
E.~Feireisl and H.~Petzeltov{\'a}.
\newblock On integrability up to the boundary of the weak solutions of the
  {N}avier-{S}tokes equations of compressible flow.
\newblock {\em Commun. Partial Differential Equations}, {\bf 25}(3-4):755--767,
  2000.
  \bibitem{fe}
 E.~Feireisl, V. M\' acha, \v S. Ne\v casov\' a, M. Tucsnak.
  \newblock Analysis of the adiabatic piston problem via methods of
continuum mechanics.
\newblock {\em Accepted to ANHP}


\bibitem{FN}
H.~Fujita, N.~Sauer.
\newblock On existence of weak solutions of the Navier-Stokes equations in regions with moving boundaries.
\newblock{\em J. Fac. Sci. Univ. Tokyo}, Sect. I {\bf 17} : 403--420, 1970.

\bibitem{Gruber2}
C.~Gruber and G.~P. Morriss.
\newblock A {B}oltzmann equation approach to the dynamics of the simple piston.
\newblock {\em J. Statist. Phys.}, {\bf 113}(1-2):297--333, 2003.

\bibitem{Gruber1}
C.~Gruber, S.~Pache, and A.~Lesne.
\newblock Two-time-scale relaxation towards thermal equilibrium of the
  enigmatic piston.
\newblock {\em J. Statist. Phys.}, {\bf 112}(5-6):1177--1206, 2003.

\bibitem{Gruber3}
C.~Gruber, S.~Pache, and A.~Lesne.
\newblock On the second law of thermodynamics and the piston problem.
\newblock {\em J. Statist. Phys.}, {\bf 117}(3-4):739--772, 2004.

\bibitem{KrMaNeWr}
O.~Kreml, V.~M\'acha, \v S.~Ne\v casov\'a and A.~Wr\'oblewska-Kami\'nska.
\newblock Weak solutions to the full Navier-Stokes-Fourier system with slip boundary conditions in time dependent domains.
\newblock {\em Journal de Math\'ematiques Pures et Appliqu\'ees.}
available online September 2017.



\bibitem{LAD2}
O.~A. Ladyzhenskaja.
\newblock An initial-boundary value problem for the {N}avier-{S}tokes equations
  in domains with boundary changing in time.
\newblock {\em Zap. Nau\v cn. Sem. Leningrad. Otdel. Mat. Inst. Steklov.
  (LOMI)}, {\bf 11}:97--128, 1968.

%\bibitem{LuMa}
%J. Luke\v s and J. Mal\'y.
%\newblock Measure and Integral
%\newblock Matfyzpress, Prague, 2005.

\bibitem{Lieb}
E.~H. Lieb.
\newblock Some problems in statistical mechanics that {I} would like to see
  solved.
\newblock {\em Phys. A}, {\bf 263}(1-4):491--499, 1999.
\newblock STATPHYS 20 (Paris, 1998).

\bibitem{LI4}
P.-L. Lions.
\newblock {\em Mathematical topics in fluid dynamics, Vol.2, Compressible
  models}.
\newblock Oxford Science Publication, Oxford, 1998.

\bibitem{DTT16}
D.~ Maity, T.~Takahashi, and M.~ Tucsnak.
\newblock {Analysis of a system modelling the motion of a piston in a viscous gas}.
\newblock{\em J. Math. Fluid Mech.} {\bf 19} 3: 551--579, 2017.
\bibitem{Neust1}
J.~Neustupa.
\newblock Existence of a weak solution to the {N}avier-{S}tokes equation in a
  general time-varying domain by the {R}othe method.
\newblock {\em Math. Methods Appl. Sci.}, {\bf 32}(6):653--683, 2009.

\bibitem{NeuPen1}
J.~Neustupa and P.~Penel.
\newblock A weak solvability of the {N}avier-{S}tokes equation with {N}avier's
  boundary condition around a ball striking the wall.
\newblock In {\em Advances in mathematical fluid mechanics}, pages 385--407.
  Springer, Berlin, 2010.

\bibitem{NeuPen2}
J.~Neustupa and P.~Penel.
\newblock A weak solvability of the {N}avier-{S}tokes system with {N}avier's
  boundary condition around moving and striking bodies.
\newblock {\em Recent Developments of Mathematical Fluid Mechanics},
Editors: Amann, H., Giga, Y., Kozono, H., Okamoto, H., Yamazaki, M.
\newblock To appear in Birkhauser, series: Advances in Mathematical Fluid Mechanics.

%\bibitem{OshFed}
%S.~Osher and R.~Fedkiw.
%\newblock {\em Level set methods and dynamic implicit surfaces}, volume 153 of
%  {\em Applied Mathematical Sciences}.
%\newblock Springer-Verlag, New York, 2003.

\bibitem{poul} L. Poul.
\newblock On dynamics of fluids in astrophysics.
\newblock {\em J. Evol. Equ.} 9, 37--66, 2009.

\bibitem{PRTR}
N.~V. Priezjev and S.M. Troian.
\newblock Influence of periodic wall roughness on the slip behaviour at
  liquid/solid interfaces: molecular versus continuum predictions.
\newblock {\em J. Fluid Mech.}, {\bf 554}:25--46, 2006.
\bibitem{S}
 J.~Saal.
 \newblock Maximal regularity for the Stokes system on noncylindrical space-time domains.
 \newblock{\em J. Math. Soc. Japan}, {\bf 58}, 3: 617–-641, 2006.
%\bibitem{QIWA}
%T.~Qian, X.-P. Wang, and P.~Sheng.
%\newblock Hydrodynamic slip boundary condition at chemically patterned
%  surfaces: {A} continuum deduction from molecular dynamics.
%\newblock {\em Phys. Rev. E}, {\bf 72}:022501, 2005.

%\bibitem{Shif}
%{\`E}.~G. Shifrin.
%\newblock On {O}. {A}. {L}adyzhenskaya's theorem on an initial-boundary value
%  problem for the {N}avier-{S}tokes equations in a domain with a time-varying
%  boundary.
%\newblock {\em Dokl. Akad. Nauk}, {\bf 418}(1):28--32, 2008.

\bibitem{Shel_first}
V.~V. Shelukhin.
\newblock The unique solvability of the problem of motion of a piston in
              a viscous gas.
\newblock {\em Dinamika Sploshn. Sredy}, ({\bf 31}):132--150, 1977.


\bibitem{Shel}
V.~V. Shelukhin.
\newblock Motion with a contact discontinuity in a viscous heat conducting gas.
\newblock {\em Dinamika Sploshn. Sredy}, ({\bf 57}):131--152, 1982.



\bibitem{StoCar}
Y.~Stokes and G.~Carrey.
\newblock On generalised penalty approaches for slip, free surface and related
  boundary conditions in viscous flow simulation.
\newblock {\em Inter. J. Numer. Meth. Heat Fluid Flow}, {\bf 21}:668--702,
  2011.

\bibitem{Wright1}
P.~Wright.
\newblock A simple piston problem in one dimension.
\newblock {\em Nonlinearity}, {\bf 19}(10):2365--2389, 2006.

\bibitem{Wright2}
P.~Wright.
\newblock The periodic oscillation of an adiabatic piston in two or three
  dimensions.
\newblock {\em Comm. Math. Phys.}, {\bf 275}(2):553--580, 2007.

\bibitem{Wright3}
P.~Wright.
\newblock {\em Rigorous results for the periodic oscillation of an adiabatic
  piston}.
\newblock ProQuest LLC, Ann Arbor, MI, 2007.
\newblock Thesis (Ph.D.)--New York University.

\end{thebibliography}
%\bibliographystyle{plain}
\end{document}